\documentclass[aoas,preprint]{imsart}

\usepackage{xspace}
\usepackage{relsize}
\usepackage{graphicx}
\usepackage{natbib}
\usepackage{subcaption}

\usepackage{amsmath,amsbsy,amstext,amssymb}
\usepackage{comment}
\usepackage{multicol}
\usepackage{algorithm, algorithmic}

\usepackage{definitions}

\newtheorem{prop}{Proposition}

\newtheorem{remark}{Remark}

\usepackage{hyperref}
\hypersetup{
	colorlinks=true,
	urlcolor=blue,
	citecolor=blue,
	linktoc=all,
	linkcolor=red} 

\arxiv{arXiv:1402.3093}

\begin{document}
\begin{frontmatter}

\title{Bayesian nonparametric dependent model for partially replicated data: the influence of fuel spills on species diversity\thanksref{T1}}
\runtitle{Dependent model for partially replicated data}
\thankstext{T1}{Supported by the European Research Council (ERC) through StG "N-BNP" 306406, by ANR BANDHITS and by the Australian Research Council.}

\begin{aug}
  \author{\fnms{Julyan}  \snm{Arbel},\thanksref{t2,m1}
  \ead[label=e1]{julyan.arbel@carloalberto.org}}
  \author{\fnms{Kerrie} \snm{Mengersen}\thanksref{m2}\ead[label=e2]{k.mengersen@qut.edu.au}}
  \and
  \author{\fnms{Judith}  \snm{Rousseau}\thanksref{m3}\ead[label=e3]{rousseau@ceremade.dauphine.fr}}
  \runauthor{J. Arbel et al.}
  \affiliation{Collegio Carlo Alberto\thanksmark{m1}, Queensland University of Technology\thanksmark{m2} and Universit\'e Paris-Dauphine\thanksmark{m3}}

\thankstext{t2}{Also at CREST, Paris, at the beginning of this project.}

\address{Address of Julyan Arbel\\
Collegio Carlo Alberto\\
Moncalieri, Italy\\
\printead{e1}}
  
\address{Address of Kerrie Mengersen\\
Queensland University of Technology\\
Brisbane, Australia\\
\printead{e2}}

\address{Address of Judith Rousseau\\
Universit\'e Paris-Dauphine\\
Paris, France\\
\printead{e3}}
\end{aug}

\maketitle
\begin{abstract}
We introduce a dependent Bayesian nonparametric model for the probabilistic modeling of membership of subgroups in a community based on partially replicated data. The focus here is on species-by-site data, \ie community data where observations at different sites are classified in distinct species. Our aim is to study the impact of additional covariates, for instance environmental variables, on the data structure, and in particular on the community diversity. To that purpose, we introduce dependence a priori across the covariates, and show that it improves posterior inference. We use a dependent version of the \GEMname distribution defined via the stick-breaking construction. This distribution is obtained by transforming a Gaussian process whose covariance function controls the desired dependence. The resulting posterior distribution is sampled by Markov chain Monte Carlo. We illustrate the application of our model to a soil microbial dataset acquired across a hydrocarbon contamination gradient at the site of a fuel spill in Antarctica. This method allows for inference on a number of quantities of interest in ecotoxicology, such as diversity or effective concentrations, and is broadly applicable to the general problem of communities response to environmental variables.
\end{abstract}

\end{frontmatter}

\section{Introduction}

This paper was motivated by the ecotoxicological problem of studying communities, or groups of species, observed as counts
of species at a set of sites, where the composition and distribution of
species may differ among sites, and for which the sites are indexed by a
contaminant. More specifically, the soil microbial data set we are focusing on in this paper was acquired at different sites of a fuel spill region in Antarctica. Although there is now much greater awareness of human impacts on the Antarctic, substantial challenges remain. One of these is the containment of historic buried station waste, chemical dumps and fuel spills. These wastes do not break down in such extreme environments and their spread is exacerbated by melting ice in summer. In order to develop effective containment strategies, it is important to understand the impact of these incursions on the natural environment. The data set considered here consists of soil microbial counts of operational taxonomic units, OTUs, as well as a site contaminant level measured by the total petroleum
hydrocarbon, TPH. Thus the aim is to model the probabilities of occurrence  associated with the
species at the different sites and to be able to interpret the impact of
the contaminant on the community as a whole or on a particular species. 

     This specific case study gives rise to a more general problem
that can be described as modeling the probability of membership of
subgroups of a community based on partially replicated data obtained
by observing different subsets of the subgroups at different levels of a
covariate. The problem can also be considered as the analysis of
compositional data in which the data points represent so called compositions, or
proportions, that sum to one. A typical example is the chemical
composition of rock specimens in the form of percentages of a
pre-specified number of elements \citep[see \eg][]{aitchison1982statistical,barrientos2013bayesian}. More
generally, the problem is endemic in many fields such as
biology, physics, chemistry and medicine.
Despite this, the solution to that problem remains a challenge. Common approaches are
typically based on parametric assumptions and require pre-specification
of the number of subgroups (e.g., species) in the community. In this
paper, we suggest an alternative that overcomes this drawback. The method is described in terms
of species for reasons of intuitiveness in
description, nevertheless, the approach is generally
applicable far beyond the species sampling framework.

     We propose a Bayesian nonparametric approach to both the
specific and general problems described above, using a covariate dependent
random probability measure as a prior distribution. Dependent extensions
of random probability measures, with respect to a covariate
such as time or position, have been extensively studied recently under
three broad constructions. First, a class of solutions is based on the
Chinese Restaurant process; see for instance \citet{caron2006bayesian,johnson2013bayesian}. These
are oriented towards in-line data collection and fast implementation.
Second, some approaches use completely random measures; see for example,
\citet{lijoi2013bayesian,lijoi2013dependent}. An appealing feature of this approach is
analytical tractability, which allows for more elaborate studying of the
distributional properties of the measures. Third, many strategies make
use of the stick-breaking representation, based on the line of research pioneerd by
\citet{maceachern1999dependent,maceachern2000dependent} which define dependent Dirichlet processes. See its plentiful variants which include  \citet{griffin2006order,griffin2011stick,dunson2007bayesian, dunson2008kernel, chung2009local} among others. The success of the stick-breaking constructions stems from their attractiveness from a computational point of view as well as their great flexibility in terms of full  support, which we prove for our model in Section~\ref{sec:full_support} of \supp. This is the approach that we follow here.

We define a dependent version of the \GEMname distribution (hereafter denoted \GEM), which is the distribution of the weights in a Dirichlet process, for modeling presence probabilities. Dependence is introduced via the covariance function of a Gaussian process, which allows dependent Beta random variables to be defined by inverse cumulative distribution functions transforms. The resulting model is not confined
to the estimation of diversity indices, but could also utilize the
predictive structure yielded by specific discrete nonparametric priors
to address issues such as the estimation of the number of new species
(subgroups) to be recorded from further sampling, the probability of
observing a new species at the $(n+m+1)$-th draw conditional on the
first $n$ observations, or of observing rare species, where by rare species one refers to species whose frequency is below a certain threshold \citep[see \eg][]{lijoi2007bayesian,favaro2012new}. 

The paper is organized as follows. 
In Section~\ref{sec:diversity} we describe our case study, review the ecotoxicological literature and background, and discuss diversity and effective concentration estimation. 
Section~\ref{sec:models} describes the Bayesian nonparametric model, posterior sampling and most useful properties of the model. Estimation results and ecotoxicological guidelines are given in Section~\ref{sec:applications}. A discussion on model considerations is given in Section~\ref{sec:considerations} and Section~\ref{sec:discussion} concludes this paper with a general discussion. Extended results, details of posterior computation and the proofs of our results are available in  \supp available as \citet{arbel2015supplementary}.
	
\section{Case study and ecotoxicological context\label{sec:diversity}}

\subsection{Case study and data\label{sec:data}}

As already sketched in the Introduction, our case study consists in a soil microbial data set acquired across a hydrocarbon contamination gradient at the location of a fuel spill at Australia’s Casey Station in East Antarctica ($110^{\circ}\, 32'$ E, $66^{\circ}\,  17'$ S), along a transect at 22 locations.  Microbes are classified as Operational Taxonomic Units (OTU), that we also generically refer to as species throughout the paper. OTU sequencing  were processed on genomic DNA using the \textsf{mothur} software package, see \citet{schloss2009fromTris}. We refer to \citet{siciliano2015} for a complete account on the data set acquisition. 
The total number of species recorded at least once at one site is 1,800+. All species were included in the estimation. However, we have noticed that it is possible to work with a subset of the data, consisting of those species with abundance over all measurements exceeding a given low threshold (say up to ten), without altering significantly the results. 
A crucial point for the subsequent analyses is that we order the species by \textit{decreasing overall abundance}, \ie species $j=1$ is the most numerous species in the whole data set. The variations of sampling across the sites explain why the species are not strictly ordered when considered site by site, see Figure~\ref{fig:comparison_DP_prop}. 

OTU measurements are paired with a contaminant called \TPHname  \citep[TPH, see][]{siciliano2014}, suspected to impact OTU diversity. The contamination TPH level recorded at each site ranges from 0 to 22,000  mg TPH/kg  soil. Ten sites were actually recorded as uncontaminated, \ie with TPH equal to zero. We call the microbial communities associated to these sites \textit{baseline communities}, and use them in order to define effective concentrations $EC_x$, see Section~\ref{sec:EC}. Although a continuous variable, TPH is recorded with ties that we interpret as due to measurement rounding. We jitter TPH concentrations with a random Gaussian noise (absolute value for the case TPH = 0) in order to account for measurement errors and to discriminate the ties. This noise can be incorporated in the probabilistic model. Reproducing estimation for varying values of the variance of the noise, moderate compared with the variability of TPH, have shown little to no alteration of the results.

\subsection{Ecotoxicological context\label{sec:ecotox}}
This paper focuses on an ecotoxicological case study where the goal is to predict the impact of a contaminant on an ecosystem. The common treatment of this question relies on toxicity tests, either on single species (called populations) or on multiple species (called communities). The need for appropriate modeling techniques is apparent due to data limitations, for instance in our case where data acquisition in Antarctica is extremely expensive. If single species modeling methods are now well comprehended, 
community modeling still lacks from theoretical evidence endorsement. There are two alternative community modeling approaches. On one hand, one can model single species independently and then aggregate the individual predictions into community predictions \citep[e.g.][]{ellis2011}. A drawback attached to the aggregation is the lack of appropriate uncertainty of the method, on top of which one necessarily lose crucial information by dismissing interplays across species. On the other hand, the response of the community as a whole is modeled, which generally entails the use of some univariate summaries of community responses, such as compositional dissimilarity \citep[e.g.][]{ferrier2006,ferrier2007} or rank abundance distributions \citep{foster2010}. Alternatively, the responses of multiple species can be modeled simultaneously \citep[e.g.][]{foster2010,dunstan2011,wang2012}.

Single species are commonly modeled through the probability of presence $p_j$ of each species $j$ as a function of the environmental parameters. The natural distribution for multiple species is the multinomial distribution, which  provides an intuitive framework when the sampling process consists of independent observations of a fixed number of species. Recent literature demonstrates the popularity of the multinomial distribution in ecology \citep[e.g.][]{fordyce2011,death2012multinomial,holmes2012dirichlet} and genomics \citep{bohlin2009,dunson2009nonparametric}. Our use of the \GEM distribution actually extends the multinomial distribution to cases where the number of species does not need be neither fixed nor known, \ie  where the prior is on infinite vectors of presence probabilities.

\subsection{Diversity\label{sec:div}}
Modeling presence probabilities provides a clear link to indices that describe various community properties of interest to ecologists, such as species diversity, richness,  evenness, \etc. The literature on diversity is extensive, not only in ecology  \citep{hill1973diversity,patil1982diversity,foster2010,colwell2012, death2012multinomial} but also in other areas of science, such as biology, engineering, physics, chemistry, economics, health and medicine \citep[see][]{borges1998family,havrda1967quantification,kaniadakis2005two}, and in more mathematical fields such as probability theory  \citep{donnelly1993asymptotic}. There are numerous ways to study the diversity of a population divided into groups, examples of predominant indices in ecology include the Shannon index $-\sum_j p_j\log p_j$, the Simpson index (or Gini index)  $1-\sum_j p_j^2$, on which we focus in this paper, and the Good index which generalizes both $-\sum_{j} p_j^\alpha\log^\beta p_j$, $\alpha,\beta\geq 0$ \citep{good1953population}. 

Diversity estimation, and more generally estimation of community indices based on species data, has been a statistical problem of interest for a long time. One of the reasons for that problem is simple and can be traced back to the high variability inherent to species data. For instance the most obvious estimators, hereafter referred to as \textit{empirical estimators}, which consist in plugging in empirical presence probabilities, \ie observed proportions $\hat p_{ij}$ of species $j$ at site $i$, suffer from that curse. Many treatments were proposed in the literature to account for this issue. An first approach is the field of occupancy modeling and imperfect detection, see for instance the monograph \citet{royle2008hierarchical}. We provide a concise description of imperfect detection modeling in Section~\ref{sec:meas-error} and do not pursue this direction here. Another approach, that we follow in this paper, consists in smoothing, or regularizing, empirical estimates. A Bayesian approach is a natural way to do so. Specifically, \citet{gill1979bayesian} show that using a Dirichlet prior distribution over $(p_1,\ldots,p_J)$ in the multinomial model with $J$ species greatly improves estimation over empirical counterparts. The reason for this is that using a prior prevents pathological behaviors due to outliers by smoothing the estimates. The smoothing is controlled by the Dirichlet parameter which can be conducted according to expert information. Compared to the framework of \citet{gill1979bayesian}, there is additional variability across sites in our case study. To instantiate this high variability of the empirical estimates of Simpson diversity, see their representation (dots) on Figure~\ref{fig:post_shannon}. However, we leverage this additional difficulty by borrowing of strength across the sites by following the intuition that neighboring sites should respond similarly to contaminant. The borrowing of strength is done by incorporating dependence across the sites in the prior distribution. 
In order not to impose the total number of species to be known a priori, we adopt a Bayesian nonparametric approach, hence extending the work by \citet{gill1979bayesian} from Dirichlet prior distributions to covariate-dependent Dirichlet process prior. This is also extending the model of \citet{holmes2012dirichlet} to a covariate-dependent setting with a priori unknown number of species. 
Note that this idea of using a Bayesian nonparametric approach as a smoothing technique for species data was recently adopted in the context of discovery probability, the probability of observing new species or species already observed with a given frequency. \citet{good1953population} proposed smoothed estimators popularized as Good--Turing estimators for discovery probabilities. Good--Turing estimators were shown to have a Bayesian nonparametric interpretation  \citep[see][]{lijoi2007bayesian,favaro2015rediscovering,arbel2015discovery}, which demonstrate the ability of Bayesian nonparametric methods to regularize  species data.

\subsection{Effective concentration\label{sec:EC}}

Highly relevant in terms of protecting an ecosystem, the \textit{effective concentration} at level $x$, denoted by $EC_x$, is the concentration of contaminant that causes $x$\% effect on the population relative to the baseline community \citep[e.g.][]{newman2012quantitative}. For example, the $EC_{50}$ is the median effective concentration and represents the concentration of a contaminant which induces a response halfway between the control baseline and the maximum after a specified exposure time. For single species studies, this is commonly assessed by an $x$\% increase in mortality. In applications with a multi species response as we are interested in this paper, it is the response of the community as a whole that is of interest. The $EC_x$  values are used to derive appropriate protective guidelines on contaminant concentrations, for instance in terms of waste, chemical dumps and fuel spills containment strategies. Currently, it is not clear how to best calculate $EC_x$ values using whole-community data. The $EC_x$  values can be defined in many ways depending on the specific aspects of interest to the ecological application. We illustrate the use of the Jaccard dissimilarity index, denoted by $\text{Jac}(X)$, one of the many dissimilarity variants available, as a measure of change in community composition. 
We defined the baseline community as the set of uncontaminated sites (ten sites), where TPH equals zero, see Section~\ref{sec:data}. The dissimilarity at TPH zero, denoted by $\text{Jac}_0$, is an estimate of the variability in community composition between uncontaminated sites. The $EC_x$ value is the smallest TPH value $X$  such that 
\begin{align}\label{eq:ECxdef}
 \text{Jac}(X)=1-(1-\text{Jac}_0)(1-x/100).
 \end{align}
In this way, $EC_0$, the TPH value for which there is no change relative to baseline, is obtained at $\text{Jac}(X)=\text{Jac}_0$, while $EC_{100}$ is obtained at $\text{Jac}(X)=1$, \ie for a TPH value such that the community composition becomes disjoint with the baseline. We see by Equation~\eqref{eq:ECxdef} that intermediate values are obtained by linear interpolation. 
The smallest TPH value is used so as to provide a conservative $EC_x$ estimate, since the dissimilarity curve is not guaranteed to  be monotonic. A particular feature of the model which allows us to follow this methodology is its ability to estimate the community composition between observed TPH values, since it is unlikely that the dissimilarity threshold $ \text{Jac}(X)$ sought in Equation~\eqref{eq:ECxdef} will coincide exactly with one of the measured TPH levels in the data. $95\%$ credible bands for $EC_x$ values were obtained in a similar fashion, \ie as the smallest and the largest values of, respectively, the $2.5\%$ and $97.5\%$ quantiles  of the  $EC_x$ value, again so as to provide conservative estimates. See Figure~\ref{fig:ECx_ECx} for an illustration of the method.

\section{Model\label{sec:models}}
\subsection{Data model\label{sec:sampling}}

We describe here the notations and the sampling process of covariate-dependent species-by-site count data. To each site $i=1,\ldots,I$ corresponds a covariate value $X_i\in\X$, where the space $\X$ is a subset of $\R^d$. We focus here on a single covariate, \ie $d=1$. The general case $d\geq 1$ is discussed in Section~\ref{sec:discussion}. 
 Individual observations $Y_{n,i}$ at site $i$ are indexed by $n=1,\ldots,N_i$, where $N_i$ denotes the total abundance, or number of observations. Observations $Y_{n,i}$ take on positive natural numbers values $j\in\{1,\ldots,J_i\}$ where $J_i$ denotes the number of distinct  species observed at site $i$. 
No hypothesis is made on the unknown total number of species $J=\max_i J_i$ in the community of interest, which might be infinite. 
We denote by $(\Xb,\Yb)$ the observations over all sites, where $\Xb=(X_i)_{i=1,\ldots,I}$, $\Yb =(\Y_{i}^{N_i})_{i=1,\ldots,I}$ and $\Y_{i}^{N_i}=(Y_{n,i})_{n=1,\ldots,N_i}$. The abundance of species $j$ at site $i$ is denoted by  $N_{ij}$, \ie the number of times that $Y_{n,i}=j$ with respect to index $n$. The relative abundance satisfies $\sum_{j=1}^{J_i} N_{ij} = N_{i}$. 

We model the probabilities of presence $\p=(\p(X_i))_{i=1,\ldots, I}=(p_j(X_i)_{j=1,2,\ldots})_{i=1,\ldots, I}$, where  $p_j(X_i)$ represents the probability of species $j$ under covariate $X_i$, by the following
\begin{equation}\label{eq:mixture_model}
Y_{n,i}\,\vert \,\p(X_i),X_i\simind \sum_{j=1}^\infty p_j(X_i)\delta_j,
\end{equation}
for $i=1,\ldots, I$, $n=1,\ldots, N_i$, where $\delta_j$ denotes a Dirac point mass at $j$. 

\subsection{Dependent prior distribution\label{sec:dep_prior}}

\begin{center}
\begin{figure}[ht!]
{\centering  
\includegraphics[width=.3\linewidth]{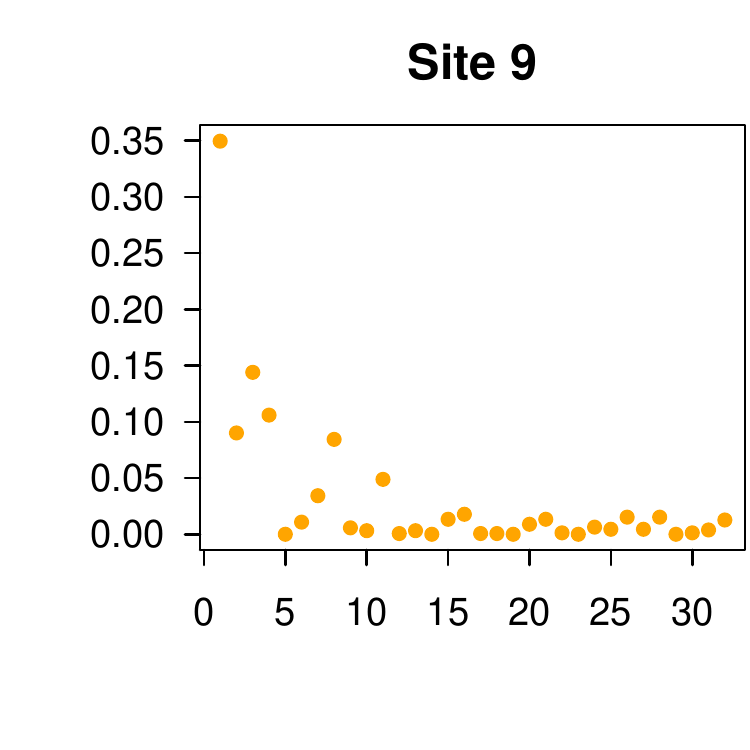} 
\includegraphics[width=.3\linewidth]{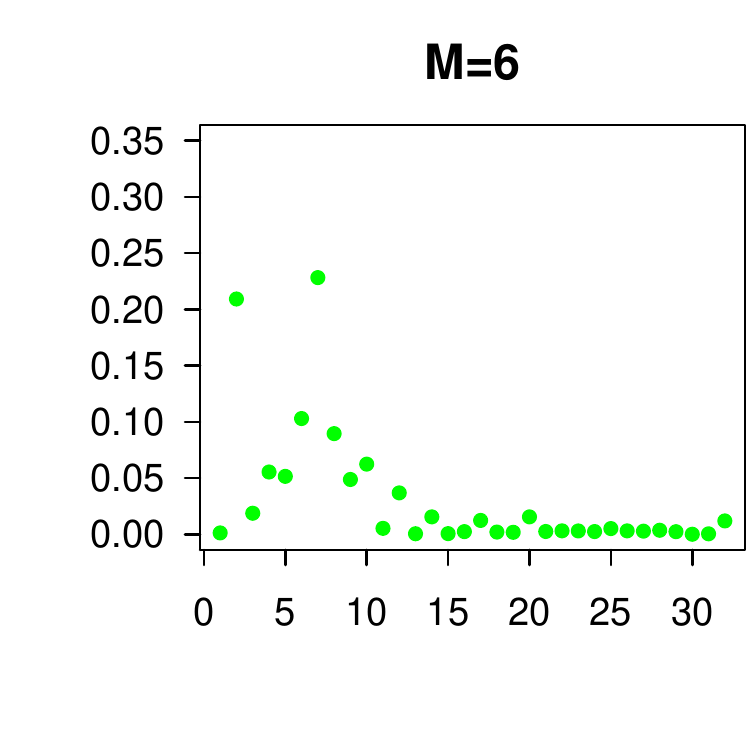} 
}
\caption[Comparison of probabilities $p_j$.]{Comparison of probabilities of presence in raw data at site $i=9$ (left) and probabilities sampled from the \gem prior with $M=6$ (right).  The $x$-axis represents  species $j=1,\ldots,32$.}
\label{fig:comparison_DP_prop}
\end{figure}
\end{center}

We follow a Bayesian approach, which implies that we need to define a prior distribution for the probabilities $\p$. The Dirichlet process \citep{ferguson1973bayesian} is a popular distribution in Bayesian nonparametrics which has been used for modeling species data by \cite{lijoi2007bayesian}. We extend the methodology developed by Lijoi et al. in building a covariate-dependent prior distribution in a way which is reminiscent of the extension of the classical \DPname to the dependent Dirichlet process by \citet{maceachern1999dependent}. More specifically, the marginal prior distribution on $\p(X)$ for covariate $X$ is defined by the following stick-breaking construction, which introduces Beta random variables $V_{j}(X)\simiid \Be(1,M)$ such that $p_1(X)=V_1(X)$ and, for $j>1$:
\begin{equation}
p_{j}(X)=V_{j}(X)\prod_{l<j}(1-V_{l}(X)).\label{eq:beta_on_V}
\end{equation}
This prior distribution is called \GEMname distribution and denoted by $\p(X)\sim\GEM(M)$, where $M>0$ is called the precision parameter. The motivation for using the \GEM distribution is explained by Figure~\ref{fig:comparison_DP_prop} which shows, for species $j=1,\ldots,32$, the observed proportions $(\hat p_{ij})$ at site $i=9$ and draws of $(p_j)$ from the $\GEM(M)$ prior with precision parameter $M=6$. 
Since the $\GEM(M)$ prior on $\p(X_i)$ is \emph{stochastically ordered} \citep[see][]{pitman2006combinatorial}, it puts more mass on the more numerous species of the community. 
It makes sense to sort the data by decreasing overall abundance, as explained in Section~\ref{sec:data}, and to use a prior with a stochastic order on $\p$ since the data under study are naturally present in large and small numbers of species. In Figure~\ref{fig:comparison_DP_prop} we observe the same non-increasing pattern between the observed frequencies and draws from the \GEM prior, which is an argument in favour of the use of the $\GEM(M)$ prior for marginal modeling of the probabilities $\p(X)$. For a discussion on the ordering assumption, see Section~\ref{sec:assump-data}.

For an exhaustive description of the prior distribution on $\p$, the marginal description~\eqref{eq:beta_on_V} needs be complemented by specifying a distribution for stochastic processes $(V_j(X),X\in\X)$, for any positive integer $j$. Since~\eqref{eq:beta_on_V} requires Beta marginals, natural candidates are Beta processes. A simple yet effective construct to obtain a Beta process is to transform a Gaussian process by the inverse cumulative distribution function (\CDF) transform as follows. Denote by $Z\sim \Norm(0,\sigma_Z^2)$ a Gaussian random variable, by $\Phi_{\sigma_{Z}}$ its \CDF  and by $F_M$ a $\Be(1,M)$ \CDF. Then $V=F_M^{-1}\circ\Phi_{\sigma_{Z}}(Z)$ is $\Be(1,M)$ distributed, with $F_M^{-1}(U) = 1-(1- U)^{1/M}$. Denote by $g_{\sigma_Z,M} =F_M^{-1}\circ\Phi_{\sigma_{Z}}$.
Note that the idea of including a transformed Gaussian process within a stick-breaking process is  used in previous articles including \citet{rodriguez2010latent,rodriguez2011nonparametric,barrientos2012support,pati2013posterior}. 

In our case, we use Gaussian processes $\Zp_j$ on the space $\X$, $j=1,2,\ldots$, which define Beta processes $\V_j$, which in turn define the probabilities $\p_j$. Though the main parameters of interest are the $\p_j$, we will work hereafter with $\Zp_j$ for computational convenience.

The Gaussian process is used as a prior probability distribution over functions. 
It is fully specified by a mean function $m$, which we take equal to 0, and a covariance function $K$ defined by
\begin{equation}\label{eq:cov_function_K}
K(X_i,X_l) = \cov\big(\Zp_j(X_i),\Zp_j(X_l)\big).
\end{equation}
We control the overall variance of $\Zp_j$ by a positive pre-factor $\sigma_{\Z}^2$ and write $K=\sigma_{\Z}^2\tilde{K}$ where $\tilde{K}$ is normalized in the sense that $\tilde{K}(X_i,X_i)=1$ for all $i$.
We work with the squared exponential (SE), Ornstein--Uhlenbeck (OU), and rational quadratic (RQ) covariance functions. See Section~\ref{sec:covariance-matrix} in \supp for more details. All three involve a parameter $\lambda$ called the length-scale of the process $\Zp_j$. It tunes how far apart two points $X_1$ and $X_2$ have to be for the process to change significantly. The shorter $\lambda$ is, the rougher are the paths of the process $\Zp_j$. 
We adopt the same technique as \citet{van2009adaptive} who deal with $\lambda$ by making it random with an inverse-Gamma (denoted \IG) prior distribution. 
They obtain adaptive minimax-optimal posterior contraction rates which indicate that the length-scale parameter $\lambda$ correctly adapts to the path smoothness. \citet{gibbs1997bayesian} derived a covariance function where the length-scale $\lambda(X)$ is a (positive) function of $X$. This case is not studied here, although it could result in interesting behaviour, as noted in \citet{Rasmussen:2006aa}.  
\begin{figure}[ht!]
\begin{center}
\includegraphics[width=.8\textwidth]{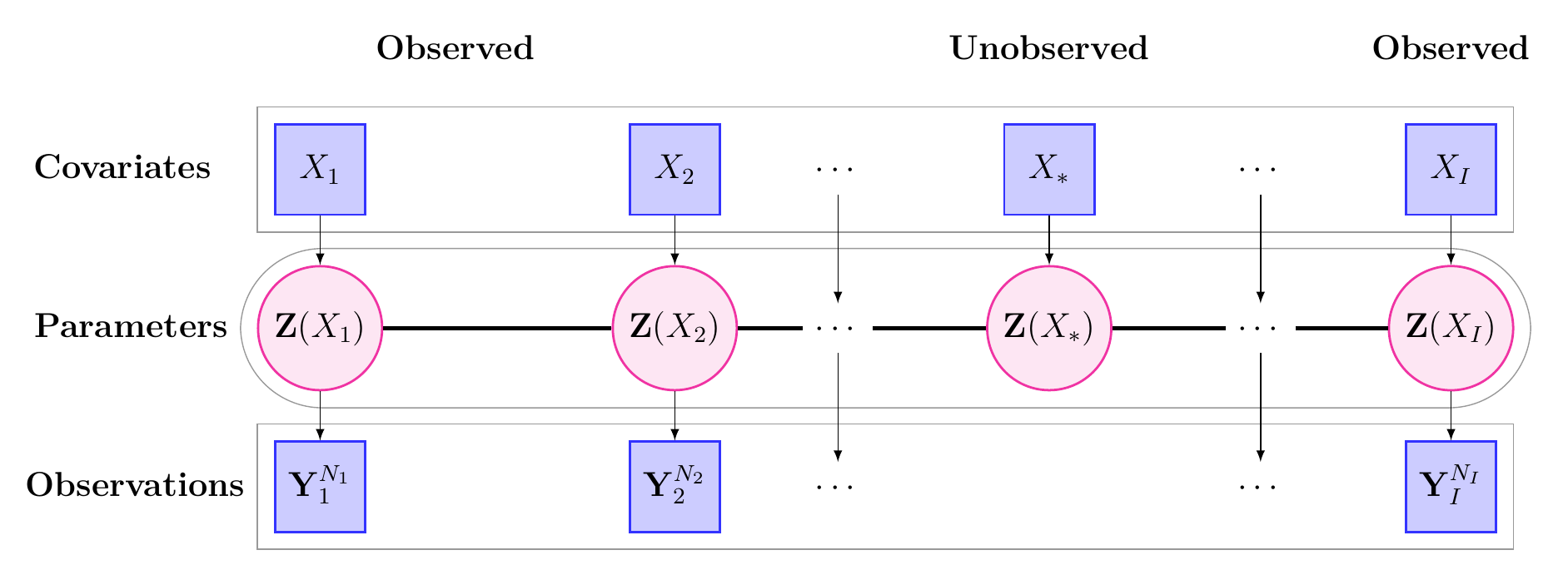}
\caption{Diagram representation for the \DGEM model. Squares represent observed data, \ie covariates $\Xb=(X_i)_{i=1,\ldots,I}$ and observations $\Y_i^{N_i}=(Y_{1,i},\ldots,Y_{N_i,i})$, and circles represent parameters for the \DGEM  model.}
\label{fig:graphical_model}
\end{center}
\end{figure}
Each species $j$ is associated to  a Gaussian process $\Zp_j$. We have a set of $I$ points $\Xb=(X_1,\ldots,X_I)$ in the covariate space $\X$ which reduces the evaluation of the whole process $\Zp_j$ to its values at $\Xb$ denoted by $\Z_j=(Z_{1,j},\ldots,Z_{I,j})= (\Zp_j(X_1),\ldots,\Zp_j(X_I))$. We denote also by $\Z$ the matrix of all vectors $\Z_j$, $\Z = (Z_{ij})_{1\leq i\leq I, 1\leq j \leq J}$. The vector $\Z_j$ is multivariate Gaussian. Its covariance matrix $K(\Xb,\lambda,\sigma_{\Z})=(\sigma_{\Z}^2\tilde{K}_{\lambda}(X_i,X_l))_{i,l=1,\ldots,I}$ is a Gram matrix with entries given by Equation~(\ref{eq:cov_function_K}). The prior distribution of $\Z_j$ is
\begin{equation*}
\log \pi(\Z_j\vert \Xb,\lambda,\sigma_{\Z}) = \frac{1}{2} \Z_j^\top K^{-1}(\Xb,\lambda,\sigma_{\Z})\Z_j - \frac{1}{2}\log\vert K(\Xb,\lambda,\sigma_{\Z})\vert - \frac{I}{2}\log 2\pi,
\end{equation*}
or, written in terms of $\sigma_{\Z}^2$ and $\tilde{K}_{\lambda}=(\tilde{K}_{\lambda}(X_i,X_l))_{i,l=1,\ldots,I}$,
\begin{align*}
\pi(\Z_j\vert \Xb,\lambda,\sigma_{\Z})\propto \sigma_{\Z}^{-I}\vert\tilde K_\lambda\vert^{-1/2}\exp\Big(-\frac{\Z_j^\top\tilde K_\lambda^{-1}\Z_j}{2\sigma_{\Z}^2}\Big).
\end{align*}
The prior distribution is complemented by specifying the distributions over  hyperparameters $\sigma_{\Z}$ the standard deviation, $\lambda$ the length-scale  and $M$ the precision parameter of the \GEM distribution. We use the following standard hyperpriors:
\begin{align}\label{eq:hyperpriors}
\sigma_{\Z}^2   \sim  \IG(a_{\Z},b_{\Z}),\,\, \lambda  \sim \IG(a_{\lambda},b_{\lambda}), \text{ and } M \sim  \Ga(a_M,b_M).
\end{align}
Note that these are also common choices in the absence of dependence since they are conjugate priors, and recall that the  inverse-Gamma for $\lambda$ also proves to lead to good convergence results. 

It is convenient to estimate the model in terms of $\Z_j$, and then to use the transform  $\V_j = g_{\sigma_{\Z},M}(\Z_j)$. The likelihood is
\begin{align}\label{eq:factorized_like}
  \Lc(\Yb\vert\Z,\Xb, \sigma_{\Z},M)=\prod_{j =1}^J\prod_{i =1}^I g_{\sigma_{\Z},M}(Z_j(X_i))^{N_{ij}} (1-g_{\sigma_{\Z},M}(Z_j(X_i)))^{\bar N_{i,j+1}},
\end{align}
where $\bar N_{i,j+1}=\sum_{l>j}{N_{il}}$. The posterior distribution is then
\begin{equation}\pi(\Z,\lambda,\sigma_{\Z},M\vert\Yb,\Xb) \propto \Lc(\Yb\vert \Z,\Xb, \sigma_{\Z},M)\pi(\Z\vert \Xb,\lambda,\sigma_{\Z})\pi(\sigma_{\Z})\pi(\lambda)\pi(M).\label{eq:post_GP}
\end{equation}
%
%

%

\subsection{Posterior computation and inference\label{sec:post}}

Here we highlight the main points of interest of the algorithm which is fairly standard, whereas the fully detailed posterior sampling procedure can be found in Supplementary Material, Section~\ref{sec:app-post}. Inference in the \DGEM model is performed via two distinct samplers: (i) first a Markov chain Monte Carlo (hereafter \MCMC) algorithm comprising Gibbs and Metropolis-Hastings steps for sampling the posterior distribution of $(\Z,\sigma_{\Z},\lambda,M)$. It proceeds by sequentially updating each parameter $\Z,\,\sigma_{\Z},\,\lambda$ and $M$ via its conditional distribution; (ii) second a sampler from the posterior predictive distribution of $\Z_*$. This consists in posterior conditional sampling of the Gaussian process $\Zp$ at covariates  $\Xb_*=(X_{1}^*,\ldots,X_{I_*}^*)$ which are not observed, \ie such that  $\{X_1,\ldots,X_I\}$ and $\{X_{1}^*,\ldots,X_{I_*}^*\}$ are pairwise distinct. This is achieved by integrating out $\Z$ in the conditional distribution of $\Z_*$ given $\Z$ according to the posterior distribution sampled in (i).

\subsection{Distributional properties\label{sec:properties}}

We provide in Proposition~\ref{prop:covariance_diversity} the first prior moments, expectation, variance and covariance, of the diversity. It is of crucial importance in order to elicit the values of hyperparameters, or their prior distribution, based on prior information (expert, etc.) Additionally, since the \DGEM introduces some dependence across the $p_j(X_i)$ in varying $X_i$, the question of the dependence induced in a diversity index arises.  Denote the Simpson index by $H_\Simp(X_i)$, see Section~\ref{sec:div}. An answer is formulated in the next Proposition in terms of the covariance between $H_\Simp(X_1)$ and $H_\Simp(X_2)$. Further properties worth mentioning are presented in \supp Section~\ref{sec:suppl-properties}, including marginal moments of the \DGEM prior and continuity of sample paths in Proposition~\ref{prop:moments}, full support in Proposition~\ref{prop:full_supp_DGEM}, a study of the joint distribution of samples from the \DGEM prior in Proposition~\ref{prop:joint_law}, and a discussion on the joint exchangeable partition probability function based on size-biased permutations in Section~\ref{sec:size-biased}.

\begin{center}
\begin{figure}[ht!]
{\centering 
\includegraphics[width=.33\linewidth]{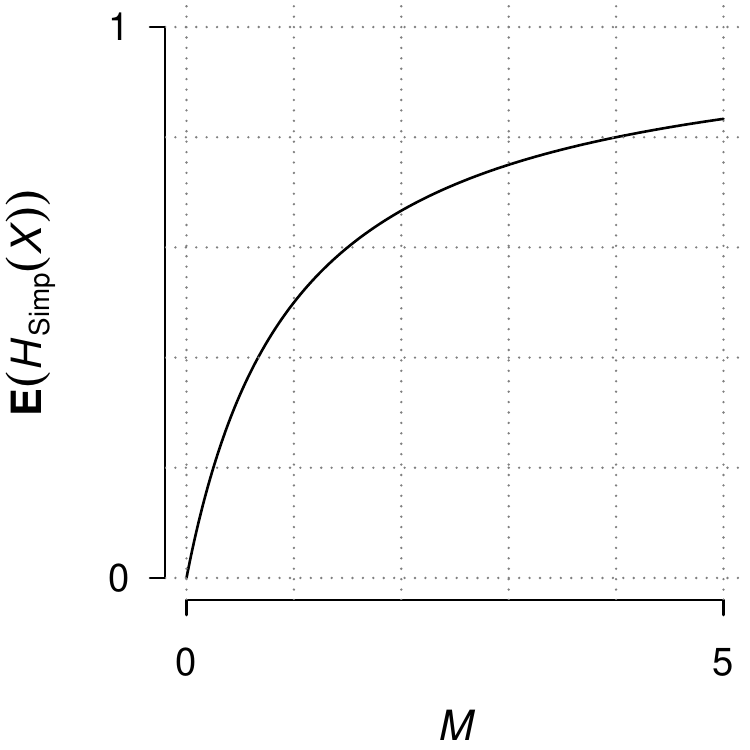}
\includegraphics[width=.33\linewidth]{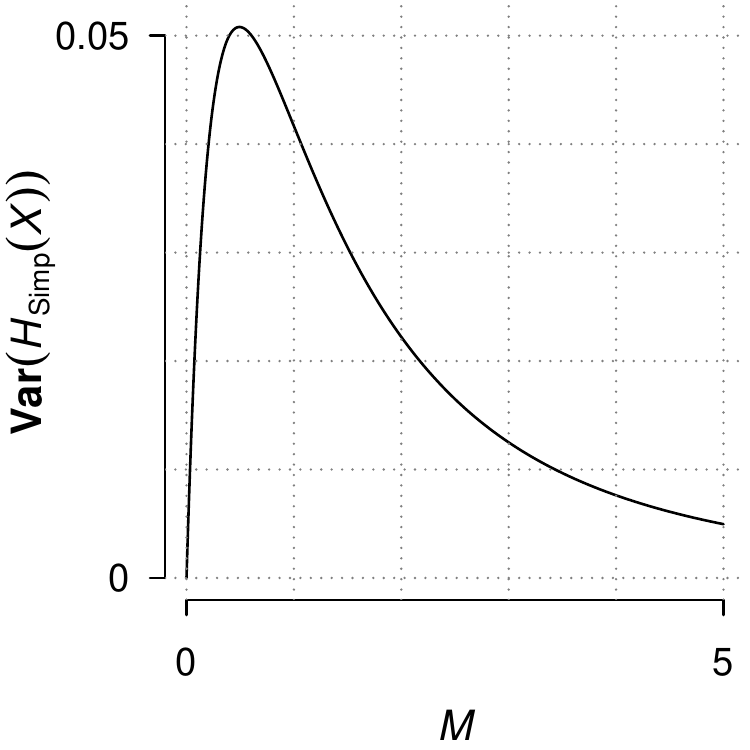}
\includegraphics[width=.32\linewidth]{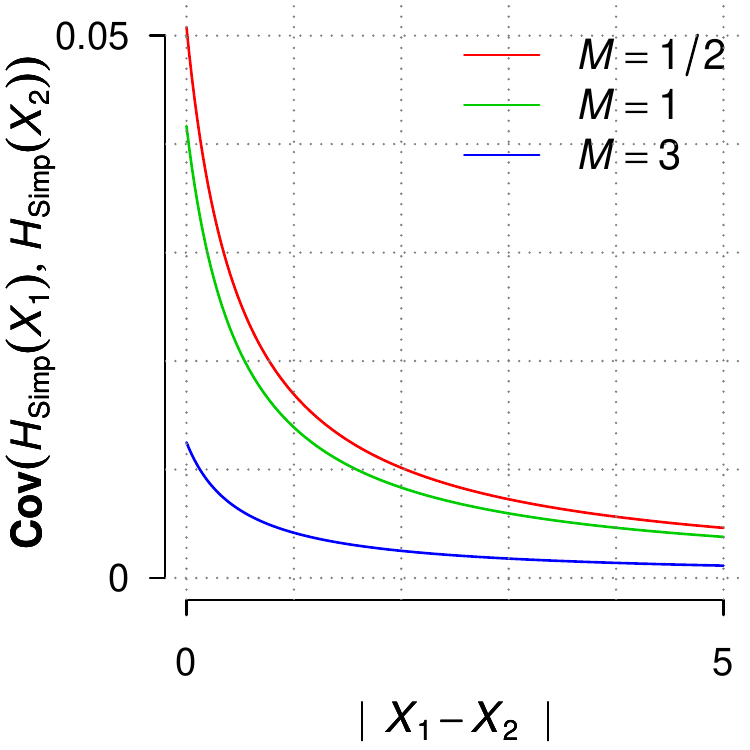}
}
\caption{Illustration of Proposition~\ref{prop:covariance_diversity}. \emph{Left}: $\E(H_\Simp(X))$ \wrt $M$. \emph{Middle}: $\var(H_\Simp(X))$ \wrt $M$. \emph{Right}: three paths of $\cov(H_\Simp(X_1),H_\Simp(X_2))$ \wrt $\vert X_1-X_2\vert$ for $M\in\MM$.}
\label{fig:var_and_cov_H}
\end{figure}
\end{center}

\begin{prop}\label{prop:covariance_diversity}
The expectation and variance of the Simpson diversity, and its covariance at two sites $X_1$ and $X_2$, induced by the \DGEM distribution, are as follows
\begin{align}
 &\E(H_\Simp)=\frac{M}{1+M},\,\var(H_\Simp(X)) \frac{2M}{(M+1)(M+1)_{3}},\label{eq:variance_simp}\\
 &\cov(H_\Simp(X_1),H_\Simp(X_2))=\frac{\nu_{2,2}(1-\omega_{2,0})+2\nu_{2,0}\gamma_{2,2}}{(1-\omega_{2,0})(1-\omega_{2,2})}-\nu_{1,0}^2,\label{eq:cov_Simpson}
\end{align}
where $\nu_{i,j}=\E[V^i(X_1)V^j(X_2)]$, $\omega_{i,j}=\E[(1-V(X_1))^i(1-V(X_2))^j]$, and $\gamma_{i,j}=\E[V^i(X_1)(1-V(X_2))^j]$.
\end{prop}
The values of $\nu_{i,j},\omega_{i,j},\gamma_{i,j}$ cannot be computed in a closed-form expression when $i\times j \neq 0$ but they can be approximated numerically. The same formal computations for the Shannon index lead to somehow more complex expressions which are not displayed here \citep[see also][]{cerquetti2014bayesian}. The expressions of Proposition~\ref{prop:covariance_diversity} are illustrated on Figure~\ref{fig:var_and_cov_H}.

The precision parameter $M$ has the following impact on the prior distribution and on the diversity: when $M\rightarrow 0$, the prior degenerates to a single species with probability 1, hence $H_\Simp\rightarrow 0$, whereas when $M\rightarrow \infty$, the prior tends to favour infinitely many species, and $H_\Simp\rightarrow 1$. In both cases, the variance and the covariance vanish. In between, the variance is maximum for $M\approx 0.49$. The covariance at $X_1$ and $X_2$ equals the variance when $X_1=X_2$ (by continuity of the sample paths), while the covariance vanishes when $\vert X_1-X_2\vert \rightarrow \infty$ (this corresponds to independence for infinitely distant covariates).
Despite the fact that the first moments of the diversity indices under a \GEM prior can be derived, a full description of the distribution seems hard to achieve. For instance, the distribution of the Simpson index involves the small-ball like probabilities $\P(\sum_j p_j^2 < a)$ for which, to the best of our knowledge, no result is known under the \GEM distribution.

\section{Case study results\label{sec:applications}}

We now apply the model to the estimation of diversity and of effective concentrations $EC_x$ as described in Section~\ref{sec:diversity}, and assess the goodness of fit of the model and its sensitivity to sampling variation.

\subsection{Results\label{sec:microbial}}

The \MCMC algorithm 
is run with squared exponential Gaussian processes for 50,000 iterations thinned by a factor of 5 with a burn-in of 10,000 iterations. The parameters of the hyperpriors~\eqref{eq:hyperpriors} are  $a_{\Z}=b_{\Z}=1$, $\eta_{\lambda}=1$, $a_{\lambda}=b_{\lambda}=1$ and  $a_M=b_M=1$. The efficiency and convergence of the \MCMC sampler was assessed by trace plots and autocorrelations of the parameters. 

\begin{figure}[ht!]
\begin{minipage}[b]{.5\linewidth}
\centering\includegraphics[width=\textwidth]{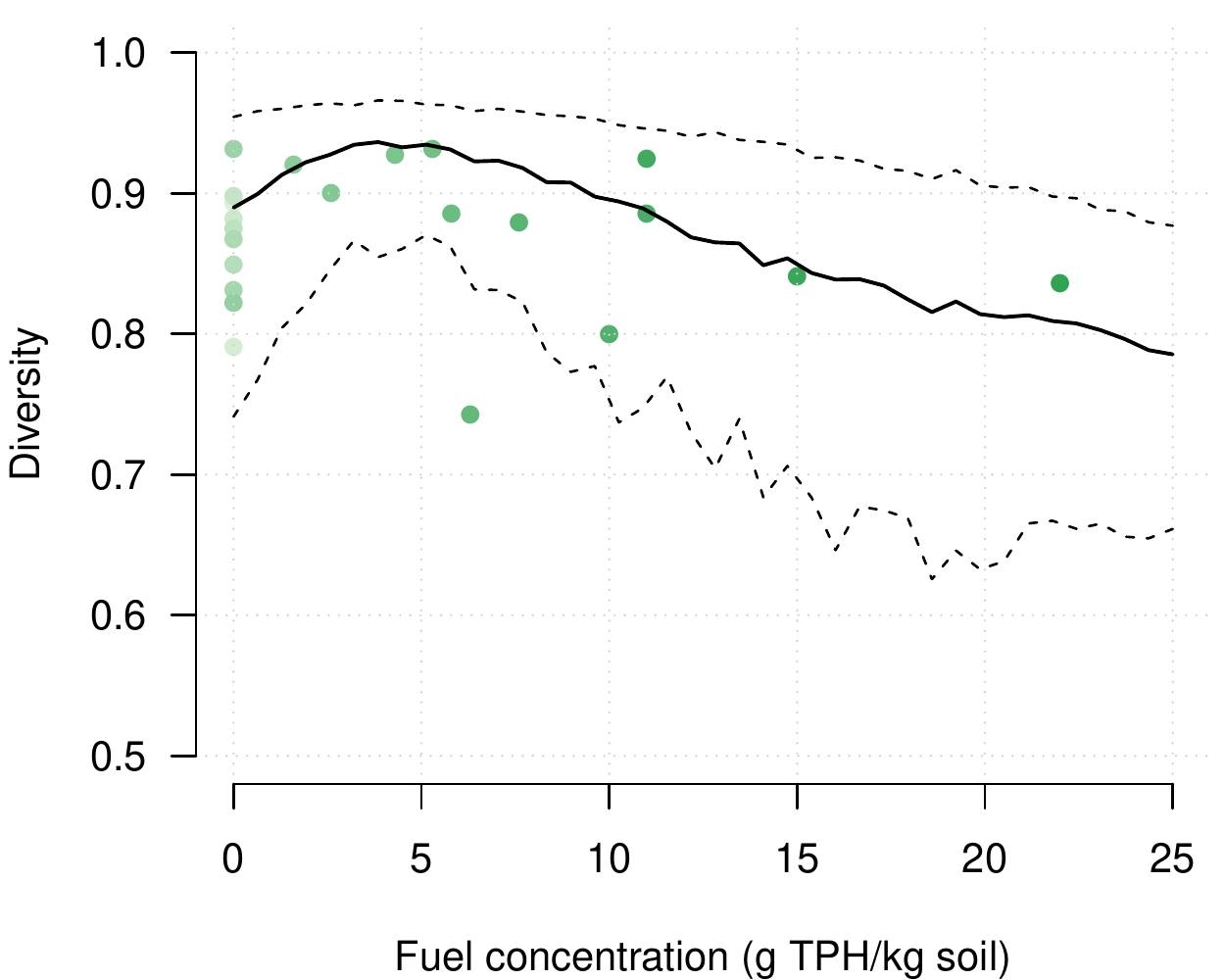}
\subcaption{\DGEM}\label{fig:diversity_DGEM}
\end{minipage}%
\begin{minipage}[b]{.5\linewidth}
\centering\includegraphics[width=\textwidth]{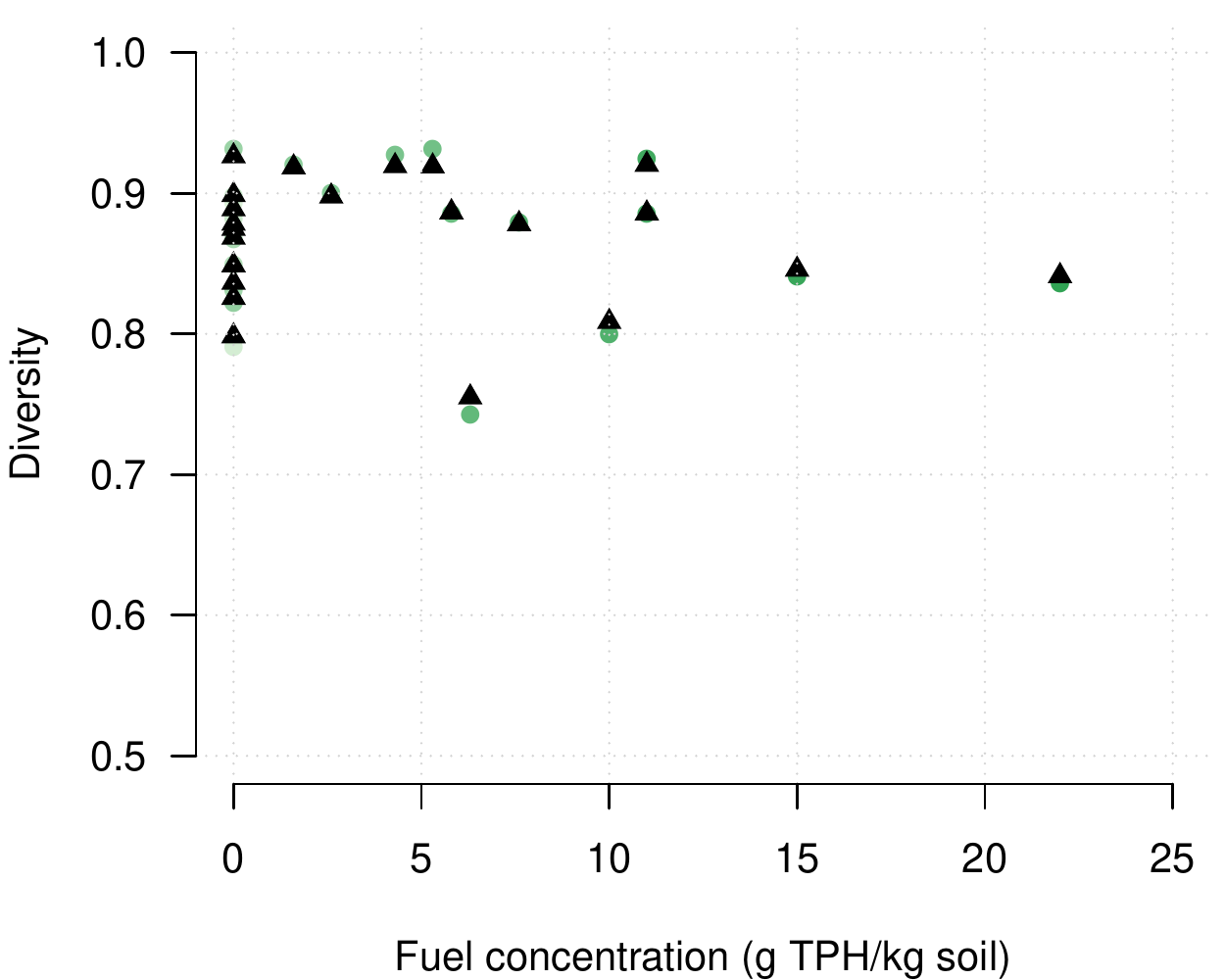}
\subcaption{Independent \GEM}\label{fig:diversity_GEM}
\end{minipage}
\caption{Diversity estimation results. {(a)} \DGEM model estimates (50,000 MCMC samples). Solid line: \SG  diversity estimate. Dashed lines: 95$\%$ credible interval for the \SG diversity.   
Dots: Empirical estimates of \SG diversity. 
{(b)} Independent \GEM model estimates  (50,000 MCMC samples). 
Triangles: posterior mean estimate of the Simpson diversity.}
\label{fig:post_shannon}
\end{figure}

The results for the Simpson diversity estimation are illustrated in Figure~\ref{fig:post_shannon} for the \DGEM model (left, \ref{fig:diversity_DGEM}) and for the independent \GEM model {(right, \ref{fig:diversity_GEM})}. The horizontal axis represents the pollution level TPH and the vertical axis represents the Simpson diversity. The posterior mean of the diversity is represented by the solid line, and a 95\% credible interval is indicated by dashed lines, for the dependent model only. The  dots indicate the empirical estimator of the diversity. 

The \DGEM model (Figure~\ref{fig:diversity_DGEM}) suggested that diversity first increases with TPH with a maximum at 4,000mg TPH/kg soil, and then decreases with TPH. The \GEM model estimates are shown for comparison in Figure~\ref{fig:diversity_GEM}. These estimates showed more variability with respect to TPH in that they are closer to the empirical estimates of the diversity. Note that the \GEM estimates were only available at levels of the covariate that were present in the data, because of the independent nature of the model specification. The \DGEM, in contrast, provided predictions across the full range of TPH values. The credible bands are narrowest for TPH between 3,000-5,000mg TPH/kg soil, due to borrowing of information between concentrated points, and they  widen both at TPH = 0, due to a lot of data points with high variability, and at large TPH, due to few data points.

\begin{figure}[ht!]
\begin{minipage}[b]{.5\linewidth}
\centering\includegraphics[width=\textwidth]{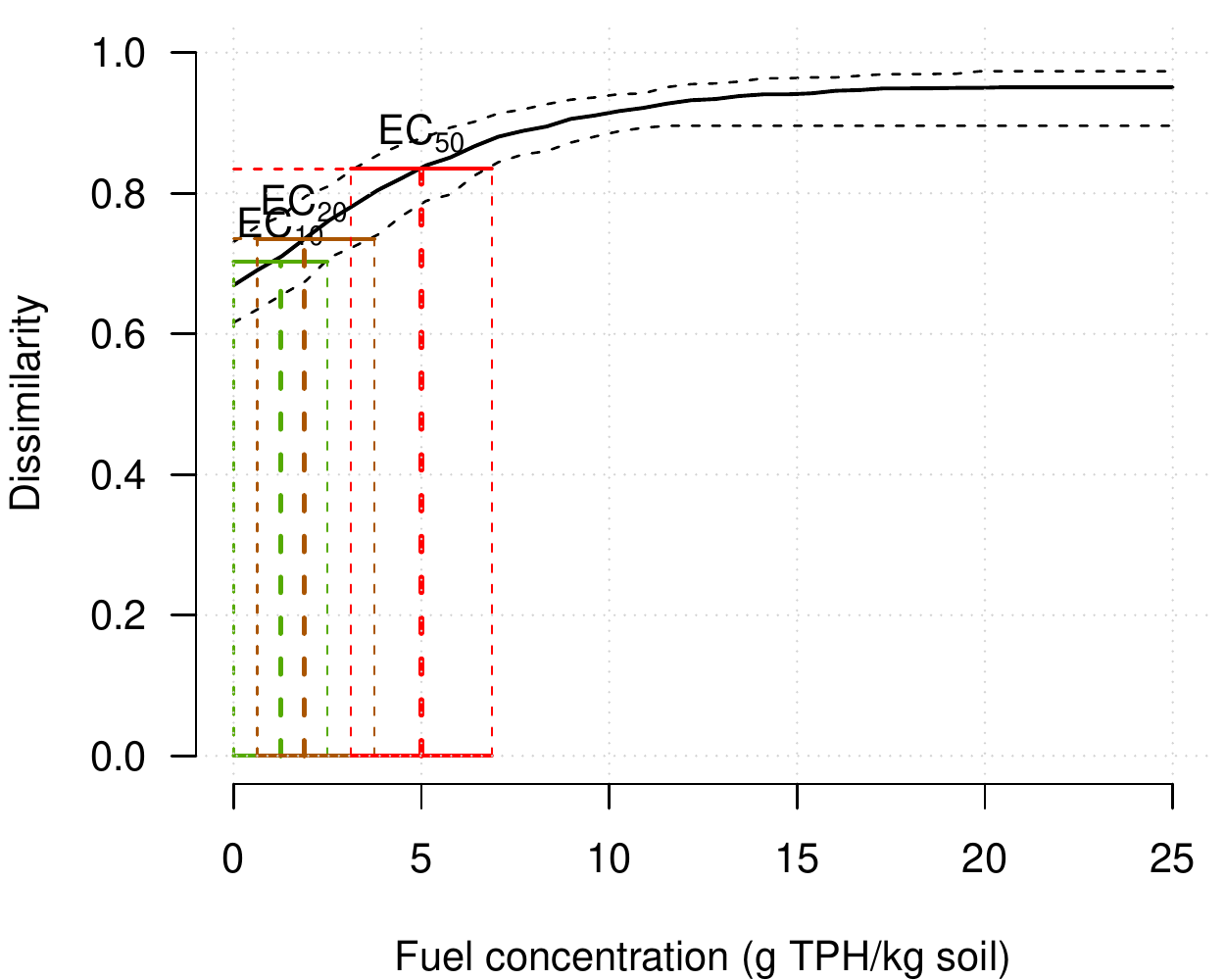}
\subcaption{Illustration of $EC_x$ and Jaccard dissimilarity}\label{fig:ECx_ECx}
\end{minipage}
\begin{minipage}[b]{.4\linewidth}
\centering
\begin{tabular}{cccc}
  \hline
$x$ & $EC_x$ & min & max \\ 
  \hline
10 & 1250 & 0 & 2500 \\ 
  20 & 1875 & 625 & 3750 \\ 
  50 & 5000 & 3125 & 6875 \\ 
   \hline
\end{tabular}

\vskip2cm
\subcaption{$EC_x$ estimates and 95\% credible intervals (min, max)}\label{tab:ECx_estimates}
\end{minipage}%
\caption{Jaccard dissimilarity and $EC_x$ estimation results. {(a)} Posterior distribution (\DGEM model) of Jaccard dissimilarity between the control community, where TPH equals zero, and communities where TPH$>0$. Solid line: mean estimate Dashed lines: 95\% credible intervals of the dissimilarity estimate. Color: Illustration of estimation of $EC_x$ values and their credible intervals. {(b)} Estimates of $EC_x$ values and their credible intervals. }
\label{fig:ECx}
\end{figure}

The Jaccard dissimilarity curve with respect to TPH is shown in Figure~\ref{fig:ECx_ECx}. The $EC_x$ values are estimated as explained in Section~\ref{sec:EC} and provided in Table~\ref{tab:ECx_estimates}. 
Dissimilarity increased with TPH, illustrating that the contaminant alters community structure. Typically, $EC_{10}$, $EC_{20}$ and $EC_{50}$ values of  Table~\ref{tab:ECx_estimates} are reported in toxicity studies to be used in the derivation of protective concentrations in environmental guidelines, see Section~\ref{sec:EC}. $EC_{10}$, $EC_{20}$ and $EC_{50}$ values estimated from this model are 1,250, 1,875 and 5,000 mg TPH/kg soil   respectively. For small $x$ (less than 10\%), the lower bound of the credible interval on the  $EC_{x}$ value is zero,  because both TPH and dissimilarity values are bounded below by zero. Conversely, for large $x$ (more than 75\%), the upper bound on the credible interval is 25,000, which is the limit of the TPH range in our analysis. 

\subsection{Posterior predictive checks\label{sec:ppc}}

\begin{figure}[ht!]
\begin{minipage}[b]{.49\linewidth}
\centering\includegraphics[width=\textwidth]{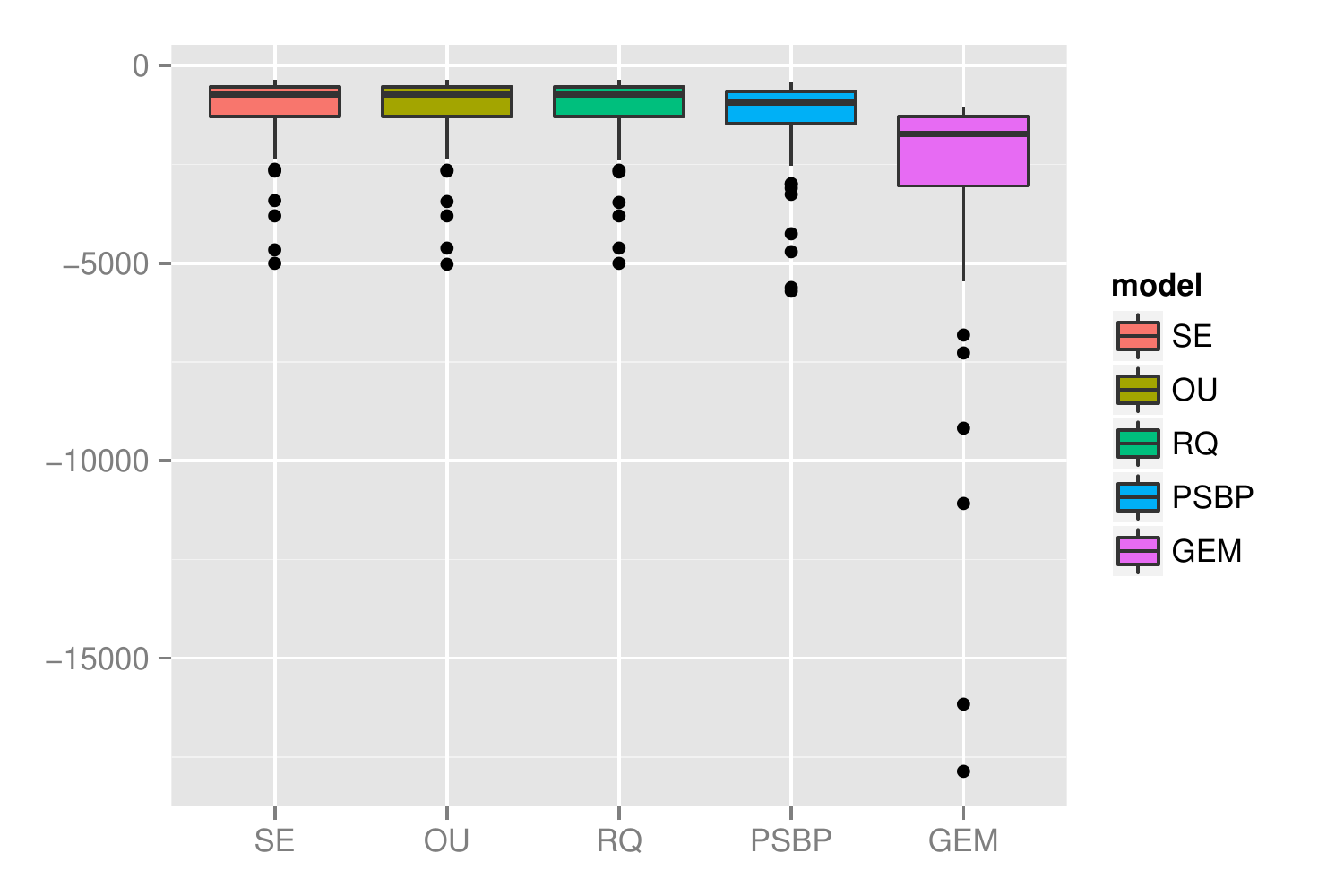}
\subcaption{}\label{fig:CPO_fig}
\end{minipage}
\begin{minipage}[b]{.49\linewidth}
\centering
\begin{tabular}{lrr}
  \hline
Model & Mean.CPO & Median.CPO \\ 
  \hline
Dep-GEM SE & -1131.3 & -732.9 \\ 
  Dep-GEM OU & -1131.5 & -732.7 \\ 
  Dep-GEM RQ & -1131.8 & -731.4 \\ 
  PSBP & -1373.5 & -932.3 \\ 
  Independent GEM & -2910.1 & -1734.1 \\ 
  Single GEM & -34606.2 & -28773.6 \\ 
   \hline
\end{tabular}

\vskip1cm
\subcaption{}\label{tab:CPO_tab}
\end{minipage}%
\caption{Log-conditional predictive ordinates (log-CPO) for different models and prior specifications (see text). {(a)} Boxplots of log-CPO. {(b)} Summaries of log-CPO, mean and median.}
\label{fig:CPO}
\end{figure}

Since we aim at comparing the performance of the model in terms of
diversity estimates, we also need to specify measures of goodness of fit. We resort to the conditional predictive
ordinates (CPOs) statistics, which are now widely used in several
contexts for model assessment. See, for example, \citet{gelfand1996model}.
For each species $j$, the CPO statistic 
is defined as follows:

\[
\operatorname{CPO}_j=\like(\Yb_j\vert \Yb_{-j}) = \int \like(\Yb_j\vert \theta) \pi(\ddr \theta\vert \Yb_{-j})
\]

where $\like$ represents the likelihood~\eqref{eq:factorized_like}, $\Yb_{-j}$ denotes data for species $j$ over all sites, $\Yb_{-j}$ denotes the  observed sample $\Yb$ with the $j$-th species excluded and $\pi(\ddr \theta \vert \Yb_{-j})$ is the posterior distribution of the model parameters
$\theta =(\Z,\sigma_{\Z},\lambda,M)$ based on data $\Yb_{-j}$. By rewriting the statistic $\operatorname{CPO}_j$ as
\[
\operatorname{CPO}_j= \biggl(\int \big(\like(\Yb_j\vert \theta )\big)^{-1} \pi(\ddr \theta \vert \Yb)
\biggr)^{-1},
\]
it can be easily approximated by Monte Carlo as
\[
\widehat{\operatorname{CPO}_j}= \Biggl(\frac{1}{T}\sum
_{t=1}^T \big(\like(\Yb_j\vert \theta^{(t)})\big)^{-1}\Biggr)^{-1},
\]
where $\{\theta^{(t)}, t=1,2,\ldots,T\}$ is an MCMC sample from
$\pi(\ddr \theta \vert \Yb)$. We illustrate the logarithm of the $\operatorname{CPO}_j$, $j=1,\ldots,J$, by boxplots in Figure~\ref{fig:CPO_fig}, and summarize their values in Table~\ref{tab:CPO_tab} in two ways, as an average of the logarithm of CPOs and as the median of the logarithm of CPOs. For the purpose of the comparison, we have estimated six models. The first three are the \DGEM model with squared-exponential (SE), Ornstein-Uhlenbeck (OU) and rational quadratic (RQ) covariance functions, see Section~\ref{sec:covariance-matrix} in \supp. The fourth is the probit stick-breaking process (PSBP) by \citet{rodriguez2011nonparametric}. For the purpose of comparison, we have set the hyperparameters of the \PSBP so as to match the expected number of clusters of the \DGEM prior. Last, we used two variants of the \GEM prior: first independent \GEM priors at each site, as in Figure~\ref{fig:diversity_GEM}, and second a single \GEM prior where the presence probabilities are all drawn from the same \GEM distribution. 

The single \GEM is used as a very crude baseline (it is not shown in the boxplots) which does poorly compared to the five other models. As expected, the dependence induced by the \DGEM and the \PSBP greatly improves the predictive quality of the model as shows the comparison to the independent \GEM. The \DGEM model has a slightly better predictive fit than the \PSBP which seems to indicate that the total ordering of the species that we use helps as far as prediction is concerned.

\subsection{Sensitivity to sampling variation\label{sec:sampling_var}}

A thorough sensitivity analysis to sampling variation was conducted in \citet{arbel2013applied}. It consisted in estimating the model on modified data, by (i) deleting the least abundant species; (ii) including additional species; (iii) excluding sites randomly. 
This sensitivity analysis showed that the model provides consistent results with data modified as described, thus supporting some robustness to sampling variation.
\section{Model considerations and extensions\label{sec:considerations}}

In addition to looking at a sensitivity analysis to sampling variation as in Section~\ref{sec:sampling_var}, here we consider sensitivity with respect to the model itself which could be extended in a number of ways.

\subsection{Imperfect detection\label{sec:meas-error}}
As pointed out in Section~\ref{sec:div} we do not connect our model to the fields of occupancy  modeling and imperfect detection developed for instance by \citet{royle2008hierarchical}. A possible extension to the current model is by accounting for imperfect detection. Following \citet{royle2006hierarchical,dorazio2008modeling}, a simple yet effective way to handle this extension is to define a probability of detection $\theta_i$ fixed for each site $i$, and to model the variability of $\theta_i$ across $i$ by an exchangeable prior. Since $\theta_i$ affects each species by the same relative proportion, the probabilities of presence $p_{j}(X_i)$ are invariant to such a formulation, and so is the diversity. Diversity being the prime focus of the present paper, we argue that there is no need to account for imperfect detection in our model, though it could be easily extended as briefly sketched if interest deviates from diversity.

\subsection{Assumption on data, stochastic decrease of the $\hat p_j$'s\label{sec:assump-data}}

We have assumed that after ordering with respect to overall abundance, the $\hat p_j$'s display a stochastically decreasing pattern as in Figure~\ref{fig:comparison_DP_prop}. In our experience, this assumption turns out to be satisfied with most of species data sets, where species can be microbes, animals, words in text, DNA sequences, \etc. However, this assumption proves to be overly restrictive in the following cases i) data might be subject to detection error: this is covered in the previous section by changing the prior adequately; ii) there are outlier species which contradict the assumption: this could be addressed by adding a mixture layer in the prior specification; iii) the underlying assumption itself is not true: this is for instance the case when all species are overall evenly distributed. A treatment would be context specific and depend on the field.

\subsection{Comparison to other models\label{sec:compar-other-models}}

In Section~\ref{sec:applications} we have compared the \DGEM model to other models: two \GEM priors and the probit stick-breaking prior (\PSBP) of \citet{rodriguez2011nonparametric}. The benefits of the \DGEM over the first two is apparent in terms of smoothing of the estimates due to the a priori dependence, see Figure~\ref{fig:post_shannon}. It also carries over better predictive fit, see Figure~\ref{fig:CPO_fig} and Table~\ref{tab:CPO_tab}, and most importantly allows us to assess the response of species to any value of the contaminant, including unsampled values. With respect to the \PSBP, the CPO indicate a slightly better predictive fit of the \DGEM prior, at least for the case study at hand.

\section{Discussion\label{sec:discussion}}

We have presented a Bayesian nonparametric dependent model for species data, based on the distribution of the weights of a \ddp, named \DGEM distribution, which is constructed thanks to Gaussian processes. A fundamental advantage of our approach based on the stick-breaking is that it brings considerable flexibility when it comes to defining the dependence structure. It is defined by the kernel of a Gaussian process, whose flexibility allows learning the different features of dependence in the data. 

In terms of model fit, we have shown that the \DGEM model improves estimation compared to an independent \GEM model. This was conducted by computing conditional predictive ordinates (CPOs). In addition, our  dependent model allows predictions at arbitrary covariate level (not just those that were in the data). It allows, for example, estimation of the diversity and the dissimilarity across the full range of covariates. This is an essential feature in applications where the experimental data are sparse and is instrumental in estimating the $EC_x$ values.

There are computational limitations to the use of this model. The estimation can deal with large number of observations since the complexity grows linearly with the number of different observed species $J$. However the number of unique covariate values $I$ represents the limiting factor of the algorithm, and may lead to dimensionality problems. 
One could consider the use of \INLA approximations \citep[see][]{rue2009approximate} in the case of prohibitively large $I$.

Possible extensions of the present paper include the following. First, extra flexibility would be guaranteed by using the two-parameter Poisson-Dirichlet distribution instead of the \GEM distribution, since it controls more effectively the posterior distribution of the number of clusters \citep{lijoi2007bayesian}. This can be done at almost no extra cost, since it only requires one additional step in the Gibbs sampler. 
Second, the \DGEM model is tested on univariate variables only, but could be extended to multivariate variables, \ie,  $X\in\R^d$, $d>1$. 
Instead of a Gaussian process $\Zp$, one would use a  Gaussian random field $\Zp^d$. To that purpose, all the methodology presented in Section~\ref{sec:models} remains valid. The algorithm can become  computationnally challenging in the case of large dimensional covariates but it does not carry additional difficulty for limited dimension. 
Applications of such an extension are promising, such as testing  joint effects in dynamical models (time $\times$ contaminant), in spatial models (position $\times$ contaminant), \etc. 

\section*{Acknowledgments}

The problem of estimating change in soil microbial diversity associated with TPH was motivated by discussions with the Terrestrial and Nearshore Ecosystems research team at the Australian Antarctic Division (AAD). The case study data used in this paper was provided by the AAD, with particular thanks to Tristrom Winsley. We acknowledge the generous technical assistance of researchers at the AAD, in particular Ben Raymond, Catherine King, Tristrom Winsley and Ian Snape. We also wish to thank Nicolas Chopin and Annalisa Cerquetti for helpful discussions, as well as the Editor, Karen Kafadar, an Associate Editor and three referees for their constructive feedback. Part of the material presented here is contained in the PhD thesis \citet{arbel2013thesis} defended at the University of Paris-Dauphine in September 2013.

\begin{supplement}
\stitle{Supplementary material}\label{suppA}
\slink[url]{Completed by the typesetter}
\sdescription{
\supplcontent 
It is postponed after the References.
}
\end{supplement}

\bibliographystyle{apalike}
\bibliography{biblio_methods,biblio_applied}
\newpage

\setcounter{page}{1}
\setcounter{section}{0}
\setcounter{equation}{0}
\renewcommand\thesection{S.\arabic{section}}
\renewcommand\theequation{S.\arabic{equation}}

\thispagestyle{empty}

\begin{center}
{\large {\textsc{\supp for \\``Bayesian nonparametric dependent model for partially replicated data: the influence of fuel spills on species diversity''}}}\bigskip

\textsc{by Julyan Arbel, Kerrie Mengersen and Judith Rousseau}\bigskip
\end{center}

\small{\supplcontent
}

\section{Posterior computation and inference in the \DGEM model\label{sec:app-post}}

Here we describe how to design a Markov chain Monte Carlo (\MCMC) algorithm  for sampling the posterior distribution of $(\Z,\sigma_{\Z},\lambda,M)$ in the \DGEM model. Up to a transformation, it is equivalent to sample the parameters in terms of Gaussian vectors $\Z$ or Beta breaks $\V$. We denote by $\pi$ the prior distribution. We make use of the factorized form of the likelihood in Equation~\eqref{eq:factorized_like} in the main paper 
in order to break the posterior sampling into $J=\max_i J_i$ independent sampling schemes. It remains a multivariate sampling scheme in terms of the $I$ sites, but avoids a very high dimensional scheme of size $I\times J$.

\subsection{\texorpdfstring{\MCMC}{MCMC} algorithm}

We use an \MCMC algorithm comprising Gibbs and Metropolis-Hastings steps for sampling the posterior distribution of $(\Z,\sigma_{\Z},\lambda,M)$, which proceeds by sequentially updating each of the parameters $\Z,\,\sigma_{\Z},\,\lambda$ and $M$ via its conditional distribution as described in Algorithm~\ref{algo:gibbs} (general sampler) and Algorithm~\ref{algo:MH} (Metropolis-Hastings step for a generic parameter $\T$). Denote by $P_{\T}(\,\cdot\,)$ the target distribution (full conditional), and by $Q_{\T}(\,\cdot\,\vertju \T)$ the proposal for a generic parameter $\T$. The variance of the latter proposal, denoted by $\sigma^2_{Q_{\T}}$, is tuned during a burn-in period.\\
\begin{minipage}[c]{5cm}
\begin{algorithm}[H]
\caption{\DGEM \label{algo:gibbs}}
\begin{itemize}
\item Update $\Z$ given $(\sigma_{\Z},\lambda,M)$
\item Update $\sigma_{\Z}$ given $(\Z,\lambda,M)$
\item Update $\lambda$ given $(\Z,\sigma_{\Z},M)$
\item Update $M$ given $(\Z,\sigma_{\Z},\lambda)$
  \end{itemize}
\end{algorithm}
\end{minipage}
\hfill
\begin{minipage}[c]{7cm}
\begin{algorithm}[H]
\caption{Metropolis-Hastings step\label{algo:MH}}
\begin{itemize}
\item  Given $\T$, propose $\T'\sim Q_{\T}(\,\cdot\,\vertju \T)$
\item Compute $\rho_{\T} = \frac{P_{\T}(\T')}{P_{\T}(\T)}\frac{Q_{\T}(\T\vert \T')}{Q_{\T}(\T'\vertju \T)}$
\item Accept $\T'$ \withproba $\min(\rho_{\T},1)$, otherwise keep $\T$
\end{itemize}
\end{algorithm}
\end{minipage}\bigskip

The full conditionals and target distributions are now fully described:
\begin{enumerate}
\item Conditional for $\Z$: Metropolis algorithm with Gaussian jumps proposal  $\Z' \sim Q_{\Z}(\,\cdot\vertju \Z) = \Norm_I(\Z,\sigma_{Q_{\Z}}^2\tilde K_\lambda)$. We use a covariance matrix proportional to the prior covariance matrix  $\tilde K_\lambda$, which leads to improved convergence of the algorithm compared to the use of a homoscedastic alternative. The target distribution is
$$P_{\Z}(\Z)\propto \Lc(\Yb\vert \Z,\Xb,\sigma_{\Z},M) \allowbreak \pi(\Z\vert \Xb,K(\Xb,\lambda,\sigma_{\Z})).$$
\item Conditional for $\sigma_{\Z}$: Metropolis-Hastings algorithm with a Gaussian proposal left truncated to 0, ${\sigma'_{\Z}} \sim Q_{{\sigma_{\Z}}}(\,\cdot\vertju {\sigma_{\Z}}) = \Norm_{\trunc}({\sigma_{\Z}},\sigma_{Q_{{\sigma_{\Z}}}}^2)$, and target distribution $$P_{\sigma_{\Z}}(\sigma_{\Z})\propto \Lc(\Yb\vert \Z,\Xb,\sigma_{\Z},M)\sigma_{\Z}^{-I-a_{\Z}/2}\exp\Big(-\frac{\Z^\top\tilde K_\lambda^{-1}\Z-2b_{\Z}}{2\sigma_{\Z}^2}\Big).$$
\item Conditional for $\lambda$: Metropolis-Hastings algorithm with a Gaussian proposal left truncated to 0, $\lambda' \sim Q_{\lambda}(\,\cdot\vertju \lambda) = \Norm_{\trunc}(\lambda,\sigma_{Q_{\lambda}}^2)$, and target distribution $$P_{\lambda}(\lambda)\propto \pi(\Z\vert \Xb,K(\Xb,\lambda,\sigma_{\Z}))\pi(\lambda).$$
\item Conditional for $M$: Metropolis algorithm with a Gaussian proposal left truncated to 0, $M' \sim Q_{M}(\,\cdot\vertju M) = \Norm_{\trunc}(M,\sigma_{Q_{M}}^2)$, and target distribution $$P_M(M)\propto M^{A_M-1}\exp(-b_M M)\prod_{i=1}^I g_{\sigma_{\Z},M}(Z_i)^{N_{ij}}(1-g_{\sigma_{\Z},M}(Z_i))^{\bar N_{i,j+1}}.$$
\end{enumerate}
\begin{remark}\label{rem:inla}
The dimensionality of the \MCMC algorithm described above equals the number of covariates $I$ (or blocks of covariates). Large dimensions can be an obstacle to the use of traditional methods (mainly due to matrix inversion). A direction that has not been investigated could be to replace \MCMC algorithms with faster approximations, of the type of \INLA for example, see \citet{rue2009approximate}.
\end{remark}

\subsection{Predictive distribution\label{sec:predictive}}

Up to now we have considered the vector $\Z$, which is the evaluation of the Gaussian process $\Zp$ at the observed covariates $\Xb=(X_1,\ldots,X_I)$. We are now interested in new outputs, called test outputs, $\Z_*$, associated with test covariates $\Xb_*=(X_{1}^*,\ldots,X_{I_*}^*)$ which are not observed, \ie  $\{X_1,\ldots,X_I\}$ and $\{X_{1}^*,\ldots,X_{I_*}^*\}$ are pairwise distinct. 
An appealing feature of the use of Gaussian processes is the possibility to easily derive the predictive distribution of $\Z_*$, which is achieved as follows. The joint distribution of the vector outputs $(\Z,\Z_*)$ according to the prior is the following $I+I_*$ multivariate Gaussian distribution
\begin{align}
\bigg(\begin{matrix}
\,\Z_{\phantom{*}}\\
\,\Z_*
\end{matrix}\bigg)
\sim
\gaussxBig{I+I_*}{
\boldsymbol{0}
}{
\bigg(\begin{matrix}
K(\Xb,\Xb) & K(\Xb,\Xb_*) \\
K(\Xb_*,\Xb) & K(\Xb_*,\Xb_*)
\end{matrix}\bigg)
},
\label{eq:gp-ads desired joint}
\end{align}
where the covariance matrices $K(\Xb,\Xb)$, $K(\Xb,\Xb_*)=K(\Xb_*,\Xb)^\top$ and $K(\Xb_*,\Xb_*)$ (resp. $I\times I$, $I\times I_*$ and $I_*\times I_*$ matrices) are defined by their entries according to the choice of the Gaussian process. 
The conditional density of $\Z_*$ given $\Z$ is the following Gaussian distribution \citep[see][]{Rasmussen:2006aa}:
\begin{align}\label{eq:cond_Z_star}
&\Z_*\vertju \Xb_*,\Xb,\Z\sim \Norm_{I_*}(m_*(\Z),K_*), \text{ with } m_*(\Z)=K(\Xb_*,\Xb) K(\Xb,\Xb)^{-1}\Z,\\
&\text{ and } K_*= K(\Xb_*,\Xb_*)- K(\Xb_*,\Xb) K(\Xb,\Xb)^{-1} K(\Xb,\Xb_*).\nonumber
\end{align}
The predictive distribution of $\Z_*$ is obtained by integrating out $\Z$ in the conditional distribution (\ref{eq:cond_Z_star}) according to the posterior distribution $\pi(\Z\vert \Yb,\Xb)$:
\begin{equation}\label{eq:pred_Z_star}
\pi(\Z_*\vertju \Xb_*,\Yb) = \int \pi(\Z_*\vertju \Xb_*,\Xb,\Z) \pi(\Z\vert \Yb,\Xb) \dd \Z.
\end{equation}
Simulating from a predictive distribution of the form of (\ref{eq:pred_Z_star}) is described in Algorithm~\ref{algo:pred}. Once a sample of $\Z$ from the posterior distribution $\pi(\Z\vert \Yb,\Xb)$ is available, one obtains a sample from the predictive distribution at almost no extra cost, by sampling from the multivariate normal distribution (\ref{eq:cond_Z_star}). 
One matrix, $K(\Xb,\Xb)$, has to be inverted, but that computation is already done for the \MCMC sampler. The variance $K_*$ of (\ref{eq:cond_Z_star}) is to be computed once. Then it is efficient to draw a sample of the desired size from the centred normal $\Norm(0,K_*)$, and then add the means $m_*(\Z)$ for $\Z$ in the posterior sample. We can obtain the predictive distribution of any $\Z_*$ associated with any test covariates $\Xb_*$, hence allowing prediction in the whole space $\X$.
\begin{center}
\begin{minipage}[c]{11cm}
  \begin{algorithm}[H]
  \caption{Predictive distribution simulation\label{algo:pred}}
\begin{itemize}
  \item  Sample $\Z$ from the posterior distribution $\pi(\Z\vert \Yb,\Xb)$
  \item  Given $\Z$, sample $\Z_*$ from the conditional distribution $\pi(\Z_*\vertju \Xb_*,\Xb,\Z)$
  \end{itemize}
  \end{algorithm}
\end{minipage}
\end{center}

\section{Covariance matrices\label{sec:covariance-matrix}}

We work with the squared exponential (SE), Ornstein--Uhlenbeck (OU), and rational quadratic (RQ) covariance functions. The next table provides the normalized covariance function $\tilde{K}(X_1,X_2)=\tilde{K}_{\lambda}(X_1,X_2)$ for these three options.
\begin{center}
\begin{tabular}{cc}
  \hline
  Covariance function & $\tilde{K}_{\lambda}(X_1,X_2)$ \\ 
  \hline
  Squared exponential (SE) & $\exp\big(- (X_1-X_2)^2/(2\lambda^2)\big)$ \\ 
   Ornstein--Uhlenbeck (OU) & $\exp\big(- \vert X_1-X_2\vert/\lambda\big)$ \\
   Rational quadratic (RQ) & $\big(1+ (X_1-X_2)^2/(2\lambda^2)\big)^{-1}$ \\  
   \hline
\end{tabular}
\end{center}

\section{Distributional properties\label{sec:suppl-properties}}

The purpose of this section is to present key distributional properties of the \DGEM prior in terms of (i) moments and continuity, (ii) full support,  (iii) dependence and (iv) size-biased permutations. Proofs are deferred to Section~\ref{sec:appendices_methodo}.

\subsection{Marginal moments and continuity}

We start by proving the continuity of sample-paths of  the process $\p\sim$\DGEM$(M)$ and providing its marginal moments. 

\begin{prop}\label{prop:moments} Let $\p\sim$\DGEM$(M)$. Then $\p$ is stationary and marginally, $\p\sim\GEM(M)$. Also, $\p$ has continuous paths (\ie $X\rightarrow (p_1(X),p_2(X),\ldots)$ is continuous for the sup norm), and its marginal moments are
\begin{align*}
&\E(p_j(X)) = \frac{M^{j-1}}{(M+1)^j},\quad \E(p_j^n(X)) = \frac{n!}{M_{(n)}}\left(\frac{M}{M+n}\right)^j,\\
&\var(p_j(X))=\frac{2M^{j-1}}{(M+1)(M+2)^{j}}-\frac{M^{2(j-1)}}{(M+1)^{2j}},\\
&\cov(p_j(X),p_k(X))=\frac{M^{(j\Max k)-1}}{(M+1)^{\vert j-k \vert +1}(M+2)^{j\Min k}}-\frac{M^{j+k-2}}{(M+1)^{j+k}},\,k\neq j,
\end{align*}
for any $j,k\geq 1$, $n\geq 0$, and where $M_{(n)}=M(M+1)\ldots (M+n-1)$ denotes the ascending factorial, $j\Max k=\max(j,k)$ and $j\Min k=\min(j,k)$.
\end{prop}

Note that the formula for $\cov(p_j(X),p_k(X))$ does not hold for $k=j$ as it does not reduce to $\var(p_j(X))$. 
The stationarity of the process as a marginal \GEM does not constrain the data to come from a stationary process. The hierarchical level of the precision parameter $M$ enables handling of diverse data structures.

\subsection{Full support of the prior\label{sec:full_support}}

The full support of the dependent Dirichlet process is proved by \citet{barrientos2012support}. Here we consider the general case of a stick-breaking prior $\Pi$ \citep{ishwaran2001gibbs} on the infinite dimensional (open) simplex
\begin{equation}\label{eq:S}
\Simplex = \left\{\p: \sum_{i=1}^\infty p_i = 1,\,\forall i \in \N^*, p_i> 0\right\}.
\end{equation}
given by $V_i\sim Be(a_i,b_i)$ iid, $a_i,b_i>0$, and $p_i=V_i\prod_{l<i}(1-V_l)$. This class of prior distributions include the \GEM distribution, as well as the distribution of the weights of the two-parameter Poisson--Dirichlet process.

\begin{prop}\label{prop:full_supp_GEM}[Full support of the \GEM prior]
For any $\epsilon>0$ and any $p^*\in \Simplex$,
$$\Pi(p:\lVert p^*-p\rVert_1 <\epsilon) >0.$$
\end{prop}

A proof can be found in \citet{bissiri2014topological}. We provide in Section~\ref{sec:appendices_methodo} another proof based on a different technique.

For the dependent GEM we introduce $ \mathcal C(\mathcal X)_+$ the set of positive and continuous functions from $\mathcal X$ to $\R$ and $\| . \|_1$ the $L_1$ norm over $\mathcal X$. 
\begin{prop}\label{prop:full_supp_DGEM}[Full support of the \DGEM prior]
Let $(V_j(X), X \in \mathcal X)$ be i.i.d stochastic processes such that  almost surely $V_j \in \mathcal C(\mathcal X)_+$, with $\mathcal X$ a compact subset of $\R^d$. Let $\mathbb P$ be the distribution of $V_j$ and $\mathbb H$ be the support of the processes $V_j$, i.e. for all $v \in \mathbb H$, 
$$ \forall \epsilon >0, \quad \mathbb P\left( \|V - v \|_1 \leq \epsilon \right) >0.$$
 Then for all $\mathbf p^\star (.) = \psi( \mathbf v^\star ) $ with $\mathbf v^\star = (v_j^\star , j \geq 1) $ and $v_j^\star \in \mathbb H $ for all $j\geq 1$ 
 $$ \pi \left( \sum_j \|  p_j - p^\star_j  \|_1  \leq \epsilon \right) >0, \quad \forall \epsilon >0$$ where $\pi$ is the distribution associated to $\mathbb P$ after the transformation $\psi$ and 
  $$  \| p_j -  p^\star_j \|_1 = \int_{\mathcal X} |p_j(x) - p_j^\star(x) | dx. $$
\end{prop}

Note that in the case  where $Z_j$ are  Gaussian processes viewed as elements of $\mathcal C([0,1])$ such as those considered in this paper, with $V_j = F_M^{-1}( \Phi_{\sigma_Z}(Z_j))$, then $ \mathbb H $ contains 
$$\left\{(p_j, j\geq 1); \, p_j \in \mathcal C([0,1]),\,  \sum_jp_j(x)=1, \, p_j(x) \geq 0 \, \forall j\geq 1\right\}.$$
  
\subsection{Joint law of a sample from the prior}

\begin{center}
\begin{figure}[ht!]
{\centering \includegraphics[width=.45\linewidth]{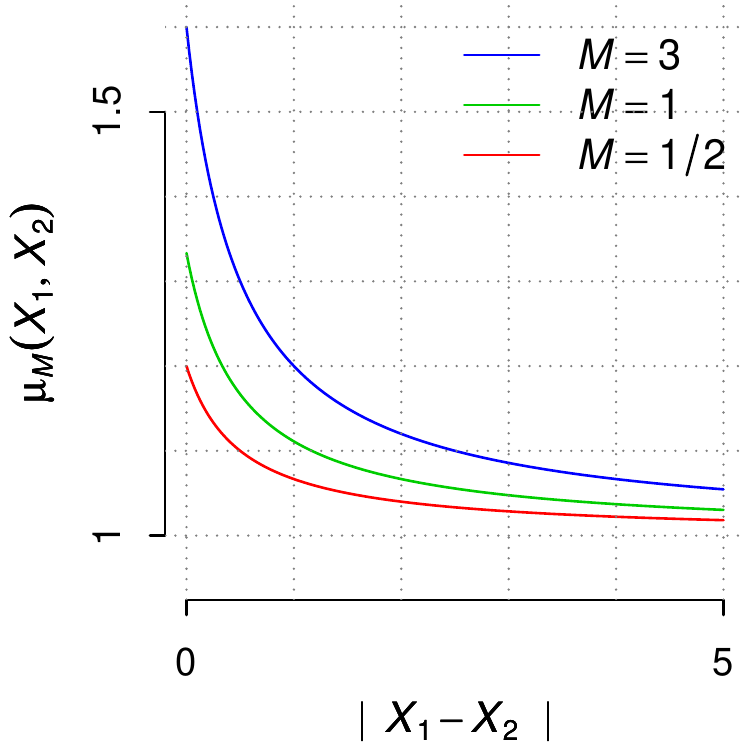} \\
}
\caption{Dependence factor $\mu_M(X_1,X_2) = (M+1)^2\E\big(V(X_1)V(X_2)\big)$ \wrt $\vert X_1-X_2\vert$ for $M\in\MM$, where $\V$ is obtained by transforming a Gaussian process with squared exponential covariance function, with $\sigma_{Z}=1$ and $\lambda=1$.}
\label{fig:mu_asymptotics}
\end{figure}
\end{center}

First, denote by $\mu_M=\mu_M(X_1,X_2)$ the dependence factor for the process evaluated at two covariates $X_1$ and $X_2$ defined by:
\begin{equation}\label{eq:mu_M}
\mu_M(X_1,X_2) = (M+1)^2\E\big(V(X_1)V(X_2)\big),
\end{equation}
Note that no analytical expression of $\mu_M$ has been derived. We resort to numerical simulation in order to compute it, \cf Figure~\ref{fig:mu_asymptotics}, and observe that $\mu_M$ is decreasing, with respect to the distance between $X_1$ and $X_2$, between two extreme cases  identified as follows:
\begin{itemize}
\item \emph{equality} case, $X_1=X_2$, \ie $V(X_1)=V(X_2)$, then $\mu_M=2(M+1)/(M+2)=1+M/(M+2)$,
\item \emph{independent} case, $V(X_1)\independent V(X_2)$ (intuitively when $\vert X_1 - X_2 \vert\rightarrow \infty$), then $\mu_M=1$.
\end{itemize} 
\begin{prop}\label{prop:joint_law}
Let observations $\Y_1^n=(Y_{1,1},\ldots,Y_{n,1})$ and $\Y_2^m=(Y_{1,2},\ldots,Y_{m,2})$ at two sites $X_1$ and $X_2$, sampled from the data model~\eqref{eq:mixture_model} conditional to the process $\p\sim $\DGEM$(M)$. The joint law of $Y_{1,1}$ and $Y_{1,2}$ is:
\begin{equation}\label{eq:joint_law_2}
\P(Y_{1,1}=j,Y_{1,2}=k) = (M+1-\mu_M)M^{\abs{j-k}-1}(M^2-1+\mu_M)^{(j\wedge k)-1}/(M+1)^{j + k},
\end{equation}
for $k\neq j$ and
\begin{equation}
\label{eq:joint_law_j_j}
\P(Y_{1,1}=j,Y_{1,2}=j) = \mu_M(M^2-1+\mu_M)^{j-1}/(M+1)^{2j},
\end{equation}
where 
$\mu_M(X_1,X_2) = (M+1)^2\E\big(V(X_1)V(X_2)\big)$ and $j\wedge k=\min(j,k)$.
\end{prop}

Equation~\eqref{eq:joint_law_2} reduces to $M^{j+k-2}/(M+1)^{j+k}$ in the \emph{independent} case (\ie $V(X_1)\independent V(X_2)$), which is indeed equal to $\P(Y_{1,1}=j)\P(Y_{1,2}=k)$.  
The probability that both first picks are equal  is obtained by summing Equation~\eqref{eq:joint_law_j_j} for all positive $j$:
\begin{align}\label{eq:joint_law_equal}
\P(Y_{1,1}=Y_{1,2}) = \frac{\mu_M}{2M+2-\mu_M}.
\end{align}
We can see that in the \emph{independent} case, Equation~(\ref{eq:joint_law_equal}) reduces to the 
probability that two draws at the same site $X_1$ belong to the same species, \ie 
$\P(Y_{1,1}=Y_{2,1}) = {1}/{(2M+1)}$, 
obtained by summing all squares of $M^{j-1}/(M+1)^{j}$.

\subsection{Size-biased permutations\label{sec:size-biased}}

In this section we derive some general results about size-biased permutations in a covariate-dependent model which are useful for the understanding of the \DGEM model. Let $\p\equ(p_1,p_2,\ldots)$ be a probability. A size-biased permutation ($\SBP$) of $\p$ is a sequence $\tilde \p=(\tilde p_1,\tilde p_2,\ldots)$ obtained by reordering $\p$ by a permutation $\sigma$ with particular probabilities. Namely, the first index appears with a probability equal to its weight, $\P(\sigma_1\equ j)= p_j$ ; the subsequent indices appear with a probability proportional to their weight in the remaining indices, \ie for $k$ distinct integers $j_1,\ldots,j_k$,
\begin{equation}\label{eq:size_biased}
\P(\sigma_k=j_{k} \vert \sigma_1= j_{1},\ldots ,\sigma_{k-1}= j_{k-1})= \frac{p_{j_k}}{1-p_{j_1}-\ldots-p_{j_{k-1}}}.
\end{equation}
We first extend Pitman's following result \citep[for example Equation~(2.23) of][]{pitman2006combinatorial}:
\begin{equation}\label{eq:pitman_lemma}
\E\Big(\sum f(p_j)\Big)=\E\Big(\sum f(\tilde p_j)\Big)=\E\bigg(\frac{f(\tilde p_1)}{\tilde p_1}\bigg),
\end{equation}
for any measurable function $f$.

\begin{prop}\label{prop:generalize_pitman}
Let $\tilde \p$ is a size-biased permutation of $\p$. For any measurable function $f$ and any integer $k\geq 1$, we have
\begin{equation}\label{eq:pitman_gene1}
\E\bigg(\sum_{(*)} f(p_{i_1},\ldots, p_{i_k})\bigg)=\E\bigg(f(\tilde p_{1},\ldots, \tilde p_{k})\prod_{i=1}^k(1-\tilde p_1-\cdots-\tilde p_{i-1})/\tilde p_i\bigg),
\end{equation}
where the sum $(*)$ runs over all distinct $i_1,\ldots,i_k$, and with the convention that the product in the right-hand side of Equation~\eqref{eq:pitman_gene1} equals $1/\tilde p_1$ when $k=1$.
\end{prop}

When it comes to averaging sums of transforms of $k$ weights $p_{i_1},\ldots, p_{i_k}$ over all distinct $i_1,\ldots,i_k$, the proposition shows that all required information is encoded by the first $k$ picks $\tilde p_{1},\ldots, \tilde p_{k}$. As stated before, the special case for $k=1$ is a well known lemma. We also mention that the case $k=2$ was proved by \citet{Archer:2013aa}.

We can look for a further insight into the \DGEM distribution by studying the \EPPFname (\EPPF) for the random variables $\Y_1^n=(Y_{1,1},\ldots,Y_{n,1})$ and $\Y_2^m=(Y_{1,2},\ldots,Y_{m,2})$ observed at covariates $X_1$ and $X_2$. See for instance \citet{pitman1995exchangeable,pitman2006combinatorial} for a summary of the importance of partition probability functions. The observations partition $[n]=\{1,2,\ldots,n\}$ and $[m]=\{1,2,\ldots,m\}$ into $k+k_1+k_2$ clusters of distinct values where
\begin{itemize}
\item $k$ clusters are commonly observed, with respective frequencies $\textbf{n}=(n_1,\ldots,n_{k})$ and $\textbf{m}=(m_1,\ldots,m_{k})$,
\item $k_1$ (resp. $k_2$) clusters are observed only at the site of covariate $X_1$ (resp. $X_2$), with  frequencies $\tilde{\textbf{n}}=(\tilde n_1,\ldots,\tilde n_{k_1})$ (resp. $\tilde{\textbf{m}}=(\tilde m_1,\ldots,\tilde m_{k_2})$).
\end{itemize}
The \EPPF can be expressed as follows
\begin{align}\label{eq:pEPPF}
p(\textbf{n},\tilde{\textbf{n}},\textbf{m},\tilde{\textbf{m}}) &=\E\bigg(\sum_{(*)} p_{i_1}^{n_1}(X_1)p_{i_1}^{m_1}(X_2)\ldots p_{i_k}^{n_k}(X_1)p_{i_k}^{m_k}(X_2) \nonumber\\
& \times p_{j_{1}}^{\tilde n_{1}}(X_1) \ldots p_{j_{k_1}}^{\tilde n_{k_1}}(X_1) \times p_{l_{1}}^{\tilde m_{1}}(X_2) \ldots p_{l_{k_2}}^{\tilde m_{k_2}}(X_2)\bigg)
\end{align}
where the sum $(*)$ runs over all  $(k+k_1+k_2)$-uples $(i_1,\ldots,i_k,j_{1},\ldots,j_{k_1},l_{1},\ldots,l_{k_2})$ with pairwise distinct elements.

In non covariate-dependent models, the \EPPF can be derived as follows. The expression of Equation~\eqref{eq:pEPPF} reduces to a simpler sum $p(\textbf{n})$ which equals the conditional expectation of the first few elements of a size-biased permutation $\tilde \p$ given $\p$, and one obtains, by application of Proposition~\ref{prop:generalize_pitman} where $f(p_{1},\ldots, p_{k}) = p_1^{n_1}\ldots p_k^{n_k}$:

\begin{equation*}
p(\textbf{n})=\E\bigg[\prod_{i=1}^k \tilde p_{i}^{n_i-1}\prod_{i=1}^{k-1}\Big(1-\sum_{j=1}^i \tilde p_j\Big)\bigg].
\end{equation*}

The \ISBPname (\ISBP) property that characterizes the \GEM distribution \citep[\cf][]{pitman1996random} can then be used to replace the first few elements of the size-biased permutation $\tilde \p$ by the first few elements of $\p$:

\begin{equation*}
p(\textbf{n})=\E\bigg[\prod_{i=1}^k p_{i}^{n_i-1}\prod_{i=1}^{k-1}\Big(1-\sum_{j=1}^i p_j\Big)\bigg].
\end{equation*}

The final steps are to use the stick-breaking representation of $\p$ with independent Beta random variables $\V$, 
and derive the \EPPF by computing the moments of Beta random variables (see Equation~\eqref{eq:moment_of_beta})
 
 \begin{equation*}
p(\textbf{n})=\frac{M^k}{M_{(n)}}\prod_{j=1}^k \left(n_j-1\right)!
\end{equation*}

Here, the hindrance to further computation of a closed-form expression for $p(\textbf{n},\tilde{\textbf{n}},\textbf{m},\tilde{\textbf{m}})$ in~\eqref{eq:pEPPF} is, to the best of our knowledge, twofold: (i) the sum in Equation~\eqref{eq:pEPPF} does not reduce to any conditional expectation of the first few elements of a size-biased permutation of $\p$, and (ii) the \ISBPname property is not straightforward to generalize to covariate-dependent distributions, hence equality in distribution between $(\tilde p_1(X_1),\tilde p_1(X_2))$ and $(p_1(X_1), p_1(X_2))$ is not a known property (whereas it is marginally true). 

Notwithstanding this, \EPPF have  been obtained in the covariate-dependent literature, though not for stick-breaking constructions, but when the dependent process is defined by normalizing random probability measures, such as completely random measures. See for instance \citet{lijoi2013bayesian,kolossiatis2013bayesian,griffin2013comparing}. See also \citet{muller2011product} for an approach based on product partition models.

\section{Proofs\label{sec:appendices_methodo}}

\subsection*{Proof of Proposition \ref{prop:moments}}

The process $\V$ constructed in the main paper 
is marginally $\Be(1,M)$, hence by the stick-breaking construction, the process $\p\sim$\DGEM$(M)$ has marginally the $\GEM(M)$ distribution. Let $\Zp\sim \GP$ as defined in the paper, 
and suppose for simplicity of notations that it is defined on $\X = \R$. Gaussian processes have continuous paths, which in turn holds for $\mathcal{V}=F_M^{-1}\circ\Phi_{\sigma_{\Zp}}(\Zp)$ since the transformation $F_M^{-1}\circ\Phi_{\sigma_{\Zp}}$ is the composition of continuous functions. Denote by $\mathcal{V}_1,\mathcal{V}_2,\ldots$ independent processes of this type, and define $\p=(\p_1,\p_2,\ldots)$ by stick-breaking,
$\p_j=\Psi_j(\mathcal{V}_1,\ldots,\mathcal{V}_j)=\mathcal{V}_j\prod_{l<j}(1-\mathcal{V}_l)$. Then for any $j$, $\Psi_j$ is a continuous function from $(0,1)^j$ to $(0,1)$, so $\p_j$ has continuous paths. This means that $\p=(\p_1,\p_2,\ldots)$  has continuous paths in the sup norm topology.

The expressions for the moments of $p_j(X)$ are derived by using the following moments of a random variable $V\sim\Be(\alpha,\beta)$, for any $j,k\geq 0$:
\begin{equation}\label{eq:moment_of_beta}
\E(V^k) = \frac{\alpha_{(k)}}{(\alpha+\beta)_{(k)}}\,\text{ and }\,
\E\big(V^k(1-V)^j\big) = \frac{\alpha_{(k)}\beta_{(j)}}{(\alpha+\beta)_{(k+j)}}.
\end{equation}
We omit the dependence in $X$ in order to simplify the notation. Note that for $V\sim \Be(1,M)$, one has $\bar V = 1-V\sim \Be(M,1)$. For any $n\geq 0$, $\E(p_j^n)$ follows from
\begin{equation}\label{eq:epjn}
\E(p_j^n)=\E\big(V_j^n\prod_{l<j}(1-V_l)^n\big)=\frac{1_{(n)}}{(M+1)_{(n)}}\Big(\frac{1_{(n)}}{(M+1)_{(n)}}\Big)^{j-1}.
\end{equation}
The formula for $\var(p_j)$ is obtained as a consequence of~\eqref{eq:epjn}, while $\cov(p_j,p_k)$, $k\neq j$, requires the computation of $\E(p_jp_k)$ as follows (suppose without loss of generality that $j>k$)
\begin{align*}
\E(V_j)&\cdot \prod_{j>l>k}\E(\bar V_l)\cdot \E(V_k \bar V_k)\cdot \prod_{k>l}\E(\bar V_l^2)\\
&= \frac{1}{M+1}\Big(\frac{M}{M+1}\Big)^{j-k-1} \Big(\frac{1}{M+1}-\frac{2}{(M+1)(M+2)}\Big) \Big(\frac{M}{M+2}\Big)^{k-1}\\
&=\frac{M^{j-1}}{(M+1)^{j-k+1}(M+2)^k}.
\end{align*}\qed

\subsection*{Proof of Proposition~\ref{prop:full_supp_GEM}}
Let $\Psi:(0,1)^{\N}\rightarrow \Simplex$ be the stick-breaking transform. It has a reciprocal defined on $\Simplex$ whose coordinates are given by 
$$V_1 = p_1,\quad V_j = p_j(1-\sum_{l=1}^{j-1} p_l)^{-1}, \quad j\geq 2,$$
which are in $(0,1)$ by construction  because for all $j$, $0<p_j<1$. 

Let $\epsilon>0$ and $p^*\in \Simplex$. Denote by $V^*$ the reciprocal of $p^*$.  Let $M = \min\{m:\lVert p_{1:m}^*\rVert_1>1-\epsilon/3\}$. 
Denote by $\Psi_M$ the restriction of $\Psi$ to its first $M$ coordinates. 
We have by construction $\Psi_M(V_{1:M}^*)=p_{1:M}^*$. Since $\Psi_M$ is continuous and $\lVert p_{1:M}^*\rVert_1>1-\epsilon/3$, there exist two neighborhoods of $V_{1:M}^*$ in $(0,1)^M$, denoted by $A_\epsilon$ and  $B_\epsilon$, such that
$$\forall V_{1:M}\in A_\epsilon,\quad \lVert p_{1:M}\rVert_1>1-\epsilon/3 \text{ for } p_{1:M} = \Psi_M(V_{1:M})$$
and
$$\forall V_{1:M}\in B_\epsilon,\quad \lVert p_{1:M}^* - p_{1:M}\rVert_1 \leq \epsilon/3  \text{ for } p_{1:M} = \Psi_M(V_{1:M})$$
The intersection of $A_\epsilon$ and $B_\epsilon$ is an open set of $(0,1)^M$ which has no trivial coordinate because it contains $V_{1:M}^*$. Denote by $D = (A_\epsilon\cap B_\epsilon)\times (0,1)^{\N}$. Then for any $V\in D$, the image $p = \Psi(V)$ satisfies
$$\lVert p-p^* \rVert_1\leq \lVert p_{1:M}-p_{1:M}^* \rVert_1 +1-\lVert p_{1:M}^* \rVert_1+1-\lVert p_{1:M} \rVert_1 \leq \epsilon$$
In addition, $D$ has positive prior mass, which proves the proposition.\qed

\subsection*{Proof of Proposition~\ref{prop:full_supp_DGEM}}

The proof follows the same line as that of Proposition~\ref{prop:full_supp_GEM}. For the sake of simplicity and without loss of generality we assume that $\int_{\mathcal X}dx = 1$.
Let $\mathbf p^\star (.) = \psi( \mathbf v^\star ) $ with $\mathbf v^\star = (v_j^\star , j \geq 1) $ and $v_j^\star \in \mathbf H $ for all $j\geq 1$. Then 
since $F_M(x) = \sum_{j=1}^M p_j^\star(x) $ is an increasing sequence (in $M$) to the constant function $1$, $\int_{\mathcal X} F_M(x) dx \uparrow 1$ and there exists $M^\epsilon $ such that 
 $$\int_{\mathcal X}F_{M^\epsilon}(x) dx \geq 1 -\epsilon/3.$$  
 The operator $\mathbf \psi_{M} : \mathbb H^M \rightarrow \mathcal C(\mathcal X)^M$ defined by 
  $ \mathbf \psi_M (V_j(.), j\leq M) = (V_j \prod_{i<j}(1-V_i)(.), j \leq M)$ is continuous for the $L_1$ norm on $\mathcal X$ for all $M$. Hence there exists an $L_1$ open neighbourhood of $(v_j^\star, j\leq M^\star)$, say $V_\epsilon$ such that  if $(v_j , j \leq M^\star) \in V_\epsilon$ 
  $$\sum_{j=1}^{M^\star} \|p_j - p_j^\star\|_1 \leq \epsilon/3, \quad (p_j, j\leq M^\star) = \mathbf \psi_{M^\star}( v_j , j \leq M^\star)$$
  the rest of the proof is the same as in the case of Proposition~\ref{prop:full_supp_GEM}. \qed
  
\subsection*{Proof of Proposition~\ref{prop:joint_law}} 

By conditional independence
\begin{align*}
\P(Y_{1,1}=j,Y_{1,2}=k) &= \E\big(\P(Y_{1,1}=j,Y_{1,2}=k\vertju \p(X_1),\p(X_2))\big)\\
&= \E(p_j(X_1)p_k(X_2)).
\end{align*}
Suppose that $j > k$, (the case $j<k$ is symmetric) then the last quantity can be decomposed into the following product of four groups of terms
\[
\begin{array}{@{}r@{}@{}c@{}@{}c@{}@{}c@{}@{}c@{}@{}c@{}@{}c@{}@{}c@{}}
&\E(V_j(X_1))&\,\cdot\,&\prod_{k<l<j}\E(\bar V_l(X_1))&\,\cdot\,& \E(\bar V_k(X_1)V_k(X_2))&\,\cdot\,&\prod_{l<k}\E(\bar V_l(X_1)\bar V_l(X_2))\\
&=\frac{1}{M+1}&\,\cdot\,&\Big(\frac{M}{M+1}\Big)^{j-k-1} &\,\cdot\,&\Big(\frac{1}{M+1}-\frac{\mu_M}{(M+1)^2}\Big) &\,\cdot\,&\Big(1-\frac{2}{M+1}+\frac{\mu_M}{(M+1)^2}\Big)^{k-1}
\end{array}
\]
which sums up to the desired quantity. 
The case $k=j$ is treated in a similar fashion.\qed

\subsection*{Proof of Proposition \ref{prop:generalize_pitman}}

By definition of the size-biased permutation, $\P(\tilde p_1 = p_i \vertju \p)=p_i$, $\P(\tilde p_2 = p_{i_2} \vertju \tilde p_1 = p_{i_1},\p)=\frac{p_{i_2}}{1-p_{i_1}}$, and  
\begin{equation}\label{eq:proba_in_pitman_generalized}
\P\big[(\tilde p_{1},\ldots, \tilde p_{k}) = (p_{i_1},\ldots, p_{i_k}) \vertju \p\big]=\prod_{l=1}^k\frac{ p_{i_l}}{1- p_{i_1}-\cdots- p_{i_{l-1}}}.
\end{equation}
Hence the right-hand side term in Proposition~\ref{prop:generalize_pitman} can be computed by double expectation and conditioning on $\p$
\begin{equation*}
\E\big[(\E\big(f(\tilde p_{1},\ldots, \tilde p_{k})\prod_{i=1}^k(1-\tilde p_1-\cdots-\tilde p_{i-1})/\tilde p_i\vertju \p \big)\big],
\end{equation*}
and a simplification arises with the probability of \eqref{eq:proba_in_pitman_generalized} when enumerating over all distinct indices  $i_1,\ldots,i_k$.\qed

\subsection*{Proof of Proposition~\ref{prop:covariance_diversity} of the main document} 

Let $\bar H(X) = 1 - H_\Simp(X) = \sum_j p_j^2(X)$. Then $\cov(H_\Simp(X_1),H_\Simp(X_2)) = \cov(\bar H(X_1),\bar H(X_2))$. First note that $\E(\bar H(X)=\E(p_1(X))=\E(V_1(X))=1/(M+1)$ by virtue of Equation~\eqref{eq:pitman_lemma}. Then $\E(\bar H(X_1)\bar H(X_2))$ is obtained by summing the following terms: 
\begin{align*}
&\text{for all }\,j\geq 1,\,\E(p_j(X_1)p_j(X_2))=\nu_{2,2}\omega_{2,2}^{j-1}, \\
&\text{for all }\,j\neq k\geq 1,\,\E(p_j(X_1)p_k(X_2))=\nu_{2,0}\omega_{2,0}^{\vert j-k\vert-1}\gamma_{2,2}\omega_{2,2}^{(j\Min k)-1},
\end{align*}
where the same kind of development as in the proof of Proposition~\ref{prop:joint_law} is employed. 
For the variance of the Simpson index, one needs, by omitting the covariate $X$ in the notation
\begin{align*}
\E\Big(\big(\sum_j p_j^2\big)^2\Big)&=\E\Big(\sum_{i,j} p_i^2 p_j^2\Big)
=\E\Big(\sum_{i\neq j} p_i^2 p_j^2\Big)+\E\Big(\sum_{i} p_j^4\Big)\\ 
&=\E(p_1(1-p_1) p_2)+\E(p_1^3)
= \E(V_1(1-V_1)^2)\E(V_2)+\E(V_1^3)\\
&=(M+6)/(M+1)_{(3)},
\end{align*}
by Proposition~\ref{prop:generalize_pitman} and the moments~\eqref{eq:moment_of_beta}.\qed

\end{document}